\documentclass[12pt]{article}

\usepackage{amsmath}
\usepackage{amssymb}
\usepackage{graphicx}

\usepackage{cite}
\usepackage[british]{babel}

\allowdisplaybreaks
  
\newcommand{\ud}{\mathrm{d}}
\newtheorem{theorem}{Theorem}[section]
\newtheorem{proposition}{Proposition}[section]
\newtheorem{lemma}{Lemma}[section]
\newtheorem{conjecture}{Conjecture}[section]
\newcommand{\tod}{\stackrel{d}{\longrightarrow}}
\newcommand{\eqd}{\stackrel{{\rm d}}{=}}

\newcommand{\inL}{\stackrel{L^1}{\longrightarrow}}
\newcommand{\inLL}{\stackrel{L^2}{\longrightarrow}}

\newcommand{\rems}{\noindent \textbf{Remarks. }}
\newcommand{\proof}{\noindent \textbf{Proof. }}
\def\Exp{{\mathbb{E}}}
\def\Pr{{\mathbb{P}}}
\def\Var{{\mathbb{V}\mathrm{ar}}}
\def\R{\mathbb{R}}

\def\1{{\bf 1 }}
\def\N{\mathbb{N}}
\def\NN{\mathcal N}
\def\X{\mathcal X}
\def\F{\mathcal F}
\def\U{\mathcal U}
\def\O{\mathcal O}
\def\tO{\tilde {\mathcal O}}
\def\Po{\mathcal P}
\def\card{{\rm card}}
\def\diam{{\rm diam}}
\def\eps{{\varepsilon}}
\def\onng{{\rm ONG}}
\def\bx{{\bf x}}
\def\by{{\bf y}}
\def\bz{{\bf z}}
\def\bX{{\bf X}}
\def\bU{{\bf U}}
\def\0{{\bf 0}}
\def\e{{\rm e}}

\def\qed{\square}

\author{Andrew R. Wade\footnote{e-mail: \texttt{Andrew.Wade@bris.ac.uk}}
\\
\normalsize
 Department of Mathematics,
 University of Bristol,\\
\normalsize
 University Walk, Bristol BS8 1TW, England.}

  \title{Asymptotic theory for the multidimensional
random on-line nearest-neighbour graph}
 \date{September 2008}

\begin{document}

\maketitle

\begin{abstract}
The on-line nearest-neighbour graph on a sequence of 
$n$
uniform random points in $(0,1)^d$ ($d \in \N$)
joins
each point after the first 
to its nearest neighbour amongst 
its predecessors.
For the total power-weighted edge-length of this graph, 
with weight exponent $\alpha \in (0,d/2]$,
we
prove $O(\max \{ n^{1-(2\alpha/d)} , \log n \})$
upper bounds on the variance. 
On the other hand, we give an $n \to \infty$
 large-sample
convergence result 
for the total power-weighted
edge-length when
 $\alpha > d/2$.  We prove corresponding
  results when the underlying point set is a Poisson
 process of intensity $n$.
\end{abstract}

\vskip 2mm

\noindent
{\em Key words and phrases:} Random spatial graphs;
 network evolution; variance asymptotics;
  martingale differences.

\vskip 2mm

\noindent
{\em AMS 2000 Mathematics Subject Classification:}  60D05  (Primary) 60F25, 90B15, 05C80 (Secondary)

\section{Introduction}
\label{int}

The (random)
on-line nearest-neighbour
graph,
which we describe in detail below,
is one of the simplest models of the evolution of (random) spatial
networks.
Graphs
with an
`on-line'
 construction, whereby vertices are added
one by one
and connected to existing vertices according to some rule,
have recently been the subject of considerable study in relation to the
modelling of real-world networks. Examples of modelling applications
include the internet, social networks, and communications networks in general. 
The literature is
extensive (see e.g.~\cite{dorog,newman} for surveys), but
mostly non-rigorous; 
rigorous mathematical results are fewer in number, even for simple
models, and existing results concentrate on graph-theoretic
rather than geometrical properties (see e.g.~\cite{bbcr,bol1}).

In recent years, much progress has been made in obtaining
large-sample
 limit theorems
for functionals defined on graphs in geometric probability, see e.g.~\cite{avbert,by,KL,penbook,mdp,mpgauss,py1,py2}. 
The graphs in question
are locally determined in a certain sense.
A natural
 functional of interest is the total (Euclidean) edge length of the graph, or, more generally,
the total power-weighted edge-length,
i.e.~the sum of the $\alpha$-powers of each edge length for a fixed 
weight exponent $\alpha>0$.
The on-line nearest-neighbour graph ($\onng$) is of particular
theoretical interest since its 
total power-weighted length functional has both normal and non-normal limiting regimes, depending on the
exponent $\alpha$. (Another example of such a graph was given in \cite{total}, but there spatial boundary
effects were crucial.) Moreover, the 
complete central limit theorem
for the ONG
seems just beyond reach of existing general results such as those of \cite{by,mdp,mpgauss,py1} which
employ various concepts of `stabilization'.

The ONG
is constructed on 
points arriving sequentially in $\R^d$ by
connecting each point (vertex)
after the first
to its nearest (in the Euclidean sense)
predecessor.
Many real-world networks have certain
characteristics in common, including spatial structure,
localization (connections tend to join nearby nodes), and sequential growth
(the network evolves over time by the addition of new nodes). The ONG is one of the simplest
models of spatial network evolution that captures these features.

The ONG appeared
in \cite{bbcr} as a  growth model
of the world wide web graph (for $d=2$), as a simplified version of the so-called
FKP network model \cite{fkp}. \cite{bbcr} studied,
amongst other things,
the vertex-degree distribution
of the ONG. Here we are concerned with geometrical properties: in particular,
the large-sample asymptotic behaviour of the
total power-weighted edge length of the $\onng$ on 
uniform
random points in the unit cube $(0,1)^d$, $d \in \N :=\{1,2,3,\ldots\}$.

In the present paper, we add to previous work on the $\onng$.
In \cite{llns}, explicit laws of large numbers were given for the total power-weighted
length of the random $\onng$ in
$(0,1)^d$, via an application
of general results from \cite{py2}.
\cite{mdp,ong} gave 
partial classification of the distributional
limits of the power-weighted length
of the $\onng$ on uniform random points in $(0,1)^d$. In particular, when $d=1$, for exponent
$\alpha>1/2$, 
\cite{ong} showed,
by a `divide-and-conquer' approach (and the `contraction method' \cite{neinrusch}),
 that the limiting
distribution of the centred total power-weighted length of the ONG
is described in terms of a
 distributional fixed-point
equation. 
In particular, these distributional limits are not Gaussian. 

It is natural to look for central limit theorems (CLTs), i.e.~proving that,
for general dimensions $d\in \N$,
for suitable values of $\alpha$, the total weight,
 centred and appropriately scaled, converges in distribution to a Gaussian limit. 
Penrose \cite{mdp} gave such a CLT
for $d \in \N$ and $\alpha \in (0,d/4)$: see Section \ref{res} below.
 As stated in \cite{mdp,ong}, it is suspected
that
 a CLT holds throughout
$\alpha \in (0,d/2]$. One contribution of the present paper is to 
give variance upper bounds for the total power-weighted
edge length of the ONG  for $\alpha \in (0,d/2]$.
These upper bounds are believed to be tight,
and are consistent with the conjectured
central limit theory. Our methods for estimating variances
are based on a martingale difference approach,
and delicate estimates of changes in the power-weighted length
of the ONG on re-sampling a particular vertex.  
 
We also give a convergence in distribution result for the total power-weighted length of the ONG,
centred as necessary, for $\alpha >d/2$. This improves on an earlier result from \cite{ong}, where
such a result was given for $\alpha >d$. We prove
this result via a refinement of the martingale
difference technique that yields the variance bounds.

Intuition behind the $\alpha = d/2$ phase transition
in the limiting behaviour is provided by the fact that increasing the weight exponent
$\alpha$ increases the relative importance of longer edges; for large enough $\alpha$
this amplifies the inhomogeneities in the structure of the ONG (`old' edges tend to be much longer)
and so destroys the Gaussian behaviour.
 
In the next section we give a formal definition of the model and state our
main results.

\section{Definitions and results}
\label{res}

Let $d \in \N$. Let  $(\bX_1, \bX_2, \ldots)$ be a sequence of
points in $(0,1)^d$.  For $n \in\N$, let $\X_n$ denote the finite
 sequence
$(\bX_1,\ldots,\bX_n)$. The on-line nearest-neighbour graph (ONG)
on vertex set $\{ \bX_1, \ldots, \bX_n\}$ is constructed by joining each point of $\X_n$
after the first by an edge
 to its nearest
 neighbour
amongst those points that precede it in the sequence. That is, for $i=2,\ldots,n$
we join $\bX_i$ by a directed edge
$(\bX_i,\bX_j)$ to $\bX_j$, $1 \leq j <i$, satisfying \[ 
\| \bX_j -\bX_i\|
= \min_{1 \leq k < i} \| \bX_k - \bX_i \|, \]
where $\|\cdot\|$ denotes the Euclidean norm on $\R^d$. We use
  lexicographic order on $\R^d$ to break any ties. 
   The resulting directed graph is the $\onng$ on $\X_n$, denoted 
$\onng (\X_n)$. 
 
 It is sometimes more convenient to view
 the ONG as an undirected graph, by
 ignoring the directedness of the edges. From
 this perspective
$\onng (\X_n)$ is a tree; in view of the
directed graph picture, it can be seen as
 rooted at $\bX_1$.

From now on we take the
points $\bX_1, \bX_2, \ldots$ to be random.
On an underlying probability space
$(\Omega,\F,\Pr)$,
let $(\bU_1,\bU_2,\ldots)$ be a sequence
of independent uniformly
distributed random vectors in $(0,1)^d$. For $n\in\N$, let 
$\U_n:=(\bU_1,\ldots,\bU_n)$. The points $\{\bU_1,\ldots,\bU_n\}$ of the 
sequence $\U_n$ then constitute a   
binomial point process
consisting of $n$ independent uniform random vectors in $(0,1)^d$.
 
For $\bx \in \R^d$ and $\X \subset \R^d$, let
$d (\bx ; \X):= \inf_{\by \in \X \setminus \{ \bx \}}
 \| \bx-\by\|$ denote
the distance from $\bx$ to its Euclidean 
nearest neighbour in $\X\setminus \{\bx\}$.
  For $d \in \N$ and $\alpha >0$, define the total
  power-weighted edge length
of 
 $\onng ( \U_n )$ by
$\O^{d,\alpha} ( \U_1 ):=0$ and for $n \geq 2$
\begin{align*}
 \O^{d,\alpha} ( \U_n ) := \sum_{i=2}^n (d(\bU_i;\U_{i-1}))^\alpha.\end{align*}
Also, define the centred version
$\tO^{d,\alpha} ( \U_n ) := \O^{d,\alpha} ( \U_n
) - \Exp [ \O^{d,\alpha} ( \U_n )
]$. We are interested in the behaviour of $\O^{d,\alpha}(\U_n)$ as $n \to \infty$.
 
We also consider the ONG defined on a Poisson number of points.
Let
$(N(t);t\geq 0)$ be the counting process of a homogeneous Poisson process
of unit rate in $(0,\infty)$, independent of $(\bU_1,\bU_2,\ldots)$.
Thus for $\lambda>0$,
$N(\lambda)$ is a Poisson random
variable with mean $\lambda$.
With $\U_n$ as defined above, for $\lambda>0$ set 
$\Po_\lambda := \U_{N(\lambda)}$. In the Poisson case, we again
use the notation
$\tO^{d,\alpha}(\Po_\lambda)=\O^{d,\alpha}(\Po_\lambda)-\Exp[ \O^{d,\alpha}(\Po_\lambda)]$
for the (deterministically) centred version.
Note that the points of the sequence $\Po_\lambda$ constitute a homogeneous (marked)
Poisson point process
of intensity $\lambda$ on $(0,1)^d$.
In this `Poissonized'
version of the ONG, we are again interested in the large-sample asymptotics,
i.e.~the limit $\lambda \to \infty$.

For $d \in \N$ let $v_d$ denote the volume of the unit-radius Euclidean
 $d$-ball, i.e.
\[ v_d :=  \pi^{d/2} \left[ \Gamma \left( 1+ (d/2) \right) \right]^{-1};\]
see e.g.~equation (6.50) of \cite{huang}.
The following result summarizes previous work
(see Theorem 4 of \cite{llns}
and Theorem 2.1 of \cite{ong}) 
on the first-order behaviour of $\O^{d,\alpha} (\U_n)$. 
Here and subsequently `$\stackrel{L^p}{\longrightarrow}$'
denotes convergence in $L^p$-norm, $p \geq 1$.

\begin{proposition} 
\label{proprop}
\cite{ong,llns}
Let $d\in\N$. For $\alpha \in (0,d)$, as $n \to \infty$
\[ n^{(\alpha-d)/d} \O^{d,\alpha} (\U_n) \inL \frac{d}{d-\alpha} v_d^{-\alpha/d} \Gamma (1+(\alpha/d)).\]
 For $\alpha=d$,  as $n \to \infty$
\[ \Exp [ \O^{d,d} (\U_n) ] \sim v_d^{-1} \log n.\]
For $\alpha > d$, there exists $\mu(d,\alpha) \in (0,\infty)$
such that as $n \to \infty$
\[ \Exp [ \O^{d,\alpha} (\U_n)] \to \mu(d,\alpha).\]
\end{proposition}
\rems
(a) In the particular case  $d=1$, Proposition 2.1 of \cite{ong} gives \[
\mu(1,\alpha)=\frac{2}{\alpha (\alpha+1)} \left( 1 + \frac{2^{-\alpha}}{\alpha-1} \right), ~~~
(\alpha > 1).\]
(b) 
These results carry over to the Poisson point process case with $\O^{d,\alpha} (\Po_n)$: this
observation
follows from 
now well-known `Poissonization' methods.\\

Second-order (i.e.~convergence in distribution) results for $\O^{d,\alpha} (\U_n)$
and $\O^{d,\alpha} (\Po_\lambda)$ were given in \cite{mdp,ong}. Specifically,
Theorem 3.6 of Penrose \cite{mdp} gives a CLT for $\alpha \in
(0,d/4)$ and Theorem 2.1(ii) of \cite{ong}
gives convergence to a non-Gaussian limit for $\alpha>d$. We summarize these results in Proposition
\ref{thm1} below. Denote by $\NN(0,\sigma^2)$ the normal distribution
with mean $0$ and variance $\sigma^2 \geq 0$; this 
includes the degenerate case $\NN(0,0) \equiv 0$. Here and subsequently
`$\tod$' denotes convergence in distribution.
\begin{proposition}
\label{thm1}
Suppose $d \in \N$. 
\begin{itemize}
\item[(i)] Suppose $\alpha \in (0,d/4)$. Then 
\cite{mdp} there exist constants
$\sigma_{d,\alpha}^2 \in [0,\infty)$ 
and $\delta_{d,\alpha}^2 \in [0,\sigma_{d,\alpha}^2]$ such that
\begin{align}
\label{1030b} \lim_{\lambda \to \infty} \lambda^{(2\alpha-d)/d} 
\Var [ \tO^{d,\alpha} (\Po_\lambda)] = \sigma_{d,\alpha}^2, ~~~
\lim_{n \to \infty} n^{(2\alpha-d)/d} 
\Var [ \tO^{d,\alpha} (\U_n)] = \sigma_{d,\alpha}^2 - \delta_{d,\alpha}^2,
\end{align}
 and as
$\lambda, n \to \infty$
\begin{align}
\label{1030a}  \lambda^{(2\alpha-d)/(2d)} \tO^{d,\alpha}(\Po_\lambda) \tod \NN (0, \sigma_{d,\alpha}^2), 
~~~ n^{(2\alpha-d)/(2d)} \tO^{d,\alpha}(\U_n) \tod \NN (0, \sigma_{d,\alpha}^2-\delta_{d,\alpha}^2)
 .\end{align} 
 \item[(ii)] Suppose $\alpha > d$. Then \cite{ong}
 there exists a mean-zero non-Gaussian random variable $Q(d,\alpha)$ such that
 as $n \to \infty$
 \begin{align}
 \label{old1}
 \tO^{d,\alpha} (\U_n) \longrightarrow Q(d,\alpha),
 \end{align}
 where the convergence is almost sure and in $L^p$, for any $p \geq 1$.
\end{itemize}
\end{proposition} 

It is conjectured (see \cite{mdp,ong}) that the
CLTs of Proposition
\ref{thm1}(i) are in fact valid
for all $\alpha \in (0,d/2)$:

\begin{conjecture}
\label{conj0}
\cite{mdp,ong}
Suppose $d \in \N$. The limit theorems
(\ref{1030b}) and (\ref{1030a})
are also valid for  $\alpha \in [d/4,d/2)$.
\end{conjecture}

In ongoing
work, we have made
some progress towards Conjecture 
\ref{conj0}, but do not yet have a proof.
 
The first main result of the present paper, Theorem
\ref{ongvar} below, provides a version of
the variance upper bounds in (\ref{1030b})  
 for all $\alpha \in (0,d/2]$. Theorem \ref{ongvar}
 is thus consistent with Conjecture \ref{conj0},
 and the bounds in Theorem \ref{ongvar} are believed
 to be sharp (up to a constant factor).
 
\begin{theorem}
\label{ongvar}
Suppose $d \in \N$.
\begin{itemize}
\item[(i)] For $\alpha \in (0,d/2)$, there is
a constant $C \in (0,\infty)$ such that for all $n \in \N$, $\lambda \geq 1$
\begin{align}
\label{a1}
 \Var [\tO^{d,\alpha} (\U_n)] \leq C n^{1-(2\alpha/d)}, ~~~ 
 \Var [\tO^{d,\alpha} (\Po_\lambda)] \leq C \lambda^{1-(2\alpha/d)}.\end{align}
\item[(ii)] 
There is
a constant $C \in (0,\infty)$ such that for all $n \in \N$, $\lambda \geq 1$
\begin{align}
\label{a3}
 \Var [\tO^{d,d/2} (\U_n)] \leq C \log (1+n), ~~~
 \Var [\tO^{d,d/2} (\Po_\lambda)] \leq C \log (1+\lambda).
\end{align}
\end{itemize}
\end{theorem}
 
 Our second main result
extends (\ref{old1}) to all $\alpha > d/2$ and also to the Poisson case.

\begin{theorem}
\label{onngthm}
Suppose $d \in \N$ and $\alpha>d/2$. Then 
there exists a mean-zero random variable $Q(d,\alpha)$
(which is non-Gaussian for $\alpha>d$) such that:
\begin{itemize}
\item[(i)] as $n\to \infty$
\begin{align}
\label{rrr}
\tO^{d,\alpha} (\U_n) \inLL Q(d,\alpha);\end{align}
\item[(ii)]
and, with the coupling of $\U_n$ and $\Po_n$ given by $\Po_n := \U_{N(n)}$,
\begin{align}
\label{333}
 \tO^{d,\alpha} (\Po_n) \inLL Q(d,\alpha).\end{align}
 \end{itemize}
\end{theorem}

\rems
(a) The fact that for $\alpha > d$ the random variables $Q(d,\alpha)$ in (\ref{old1}) and
Theorem \ref{onngthm} are
not normal follows since convergence also 
holds without any centring; 
see Theorem 2.1(ii) of \cite{ong}.
In the special case $d=1$, a weaker version of (\ref{rrr}), with convergence in distribution only,
 was given
for $\alpha>1/2$ in Theorem 2.2 of
\cite{ong}. 
In the $d=1$ case, more information can be obtained about the distribution of
 $Q(1,\alpha)$ using a `divide-and-conquer' technique; see \cite{ong},
in particular Theorem 2.2, where the distribution of $Q(1,\alpha)$, $\alpha>1/2$
 is given (in the binomial setting, and the result 
  carries over to the Poisson
 setting by Theorem \ref{onngthm} here). Indeed, $Q(1,\alpha)$, $\alpha >1/2$,
 is given by the 
 unique solution to a distributional
 fixed-point equation, and in particular is not Gaussian; 
 see \cite{ong} for details.
We suspect  
that $Q(d,\alpha)$ is non-Gaussian for $\alpha \in (d/2,d]$ also for $d \geq 2$.\\
(b) A closely related `directed' version of the one-dimensional $\onng$ is the
`directed linear tree' 
introduced in \cite{total}, in which each point in a sequence
of points in $(0,1)$
is joined
to its nearest predecessor to the {\em left}. 
Following the methods of the present
 paper, one can obtain results for that model
 analogous to the $d=1$ cases of
 all those in this section.\\

Theorem \ref{ongvar}(ii) suggests that the
 case $\alpha=d/2$ is of a special nature.
 Moreover, the case $d=2$, $\alpha=1$ is of 
natural
interest, where we have the total Euclidean length of the $\onng$
on random points in $(0,1)^2$. 
We conjecture the following.

\begin{conjecture}
\label{conj1}
Let $d \in \N$. There exists a constant $\sigma_{d,d/2}^2 \in (0,\infty)$
such that
\begin{align*}
  (\log n)^{-1/2} \tO^{d,d/2}(\U_n) 
 \tod
\NN(0,\sigma_{d,d/2}^2), ~{\rm as}~
n \to \infty. \end{align*} 
\end{conjecture}

The proof (or refutation)
of Conjecture \ref{conj1} seems to be a challenging open problem.

The structure of the remainder of the paper is as follows.
In Section \ref{prel} we give some preparatory results
on the properties of the ONG.
In Section \ref{vars} we use a martingale difference technique
to prove Theorem \ref{ongvar}. 
In Section \ref{prf} we refine the martingale difference technique to
give a proof of Theorem \ref{onngthm}.

\section{Preliminaries}
\label{prel}

 First we introduce some more notation.
Let $\card(\X)$ denote the cardinality (number of elements) of a finite
set $\X$, and let $\0$ be the origin of $\R^d$ ($d \in \N$). For measurable
$R \subset \R^d$, let $|R|$ denote the $d$-dimensional
 Lebesgue measure of $R$. Let $\diam(R)=\sup_{\bx, \by \in R}
\| \bx-\by \|$ 
denote the (Euclidean)
diameter of a bounded set $R \subset \R^d$.
Let $B(\bx;r)$ be the (closed)
Euclidean $d$-ball with centre $\bx \in\R^d$ and radius $r>0$.
 
In the analysis in Sections \ref{vars} and \ref{prf}
below, we will need
detailed properties
of the change in total weight of the ONG on $\U_n$
when the point $\bU_i$, $i \in \{1,\ldots,n\}$,
is independently re-sampled, i.e., replaced by
an independent copy $\bU_i'$. The changes due to
edges {\em incident} to $\bU_i, \bU_i'$ require
most work to deal with. To study these, we make use
of the fact that an edge from
$\bU_j$ with $j >i$
can only be incident to $\bU_i$ if $\bU_j$
falls in the Voronoi cell 
of $\bU_i$ with respect to $\{ \bU_1, \ldots, \bU_i \}$.
Hence  the preliminary results in this
section begin with an analysis of such Voronoi cells.

The next lemma
 gives bounds on the expected
diameter of Voronoi cells in $(0,1)^d$ with respect to
$\U_n$. 
 For $n \in \N$, let $V_n(\bx)$ be the Voronoi cell of $\bx \in (0,1)^d$ 
 with respect to 
$\{\bx,\bU_1,\ldots,\bU_{n}\}$:
\begin{align}
\label{voron}
 V_n(\bx) := \left\{ \by \in (0,1)^d : \|\bx-\by\| \leq \min_{1 \leq i \leq n}
\| \by -\bU_i \| \right\} \subseteq (0,1)^d.\end{align}
\begin{lemma}
\label{kkjj}
Let $d \in \N$, $\beta>0$.
Then there exists $C \in (0,\infty)$ 
such that for 
all $n \in \N$
\[ \sup_{\bx \in (0,1)^d} \Exp [ (\diam(V_n(\bx)))^\beta ] \leq C n^{-\beta/d}.\]
\end{lemma}

We will  prove Lemma \ref{kkjj} using a construction of overlapping
and nested cones from p.~1027 of \cite{mpgauss}. 
The argument works for an arbitrary convex set, not just $(0,1)^d$,
but here we only need the latter.

For $d \in \{2,3,\ldots\}$, we can (and do) choose $I \in \N$
and construct $C_i$, $1 \leq i \leq I$  a finite collection of infinite
closed cones in $\R^d$ with angular radius $\pi/12$ and apex at $\0$,
with $\cup_{i =1}^I C_i  = \R^d$. 
Let $C_i(\bx)$ be the translate of
$C_i$ with apex at $\bx \in \R^d$.
Let $C_i^+(\bx)$ be the closed
 cone with apex and principal axis coincident with
those of $C_i (\bx)$
but with angular radius $\pi/6$. When $d=1$, we take $I=2$ and let
$C_1 = [0,\infty)$, $C_2 = (-\infty,0]$, and for $\bx \in \R$
set $C_1(\bx) =C_1^+(\bx)= [\bx, \infty)$
and
$C_2(\bx) =C_2^+(\bx)= (-\infty,\bx]$.

Let $d \in \N$.
For $\bx \in \R^d$ and $r>0$, let $C_i (\bx;r) := C_i (\bx) \cap B(\bx ; r)$
and
$C_i^+ (\bx;r) := C_i^+ (\bx) \cap B(\bx ; r)$.
For $n \in \N$, define the event
\[ E_n (\bx; r) := \bigcap_{i:1 \leq i \leq I,
\diam(C_i(\bx;r) \cap (0,1)^d) = r } \{ \U_n \cap C_i^+ (\bx; r) \neq \emptyset \} ,\]
with the convention that an empty intersection is $\Omega$.
Then $E_{n} (\bx; r) \subseteq E_{n+1} (\bx ; r)$, and
for $s \geq r$, $E_n (\bx; r) \subseteq E_n (\bx;s)$. 
For $\bx \in (0,1)^d$,
 set
\begin{align}
\label{tttt}
R_n(\bx) := \inf   \{ r >0 : E_n(\bx;r) ~{\rm occurs } \} 
 .\end{align}
 Note that a.s., $R_n(\bx) \leq d^{1/2}$. The next lemma is the main step in the proof
 of Lemma \ref{kkjj}.

\begin{lemma}
\label{kkjjx}
Suppose $d \in \N$.
For $\beta>0$ 
there exists $C \in (0,\infty)$ such that for all $n \in \N$
\[ \sup_{\bx \in (0,1)^d} \Exp [ R_n(\bx) ^\beta ] \leq C n^{-\beta/d}.\]
\end{lemma}
\proof  For $\bx \in (0,1)^d$ and $r >0$,
$\Pr ( R_n (\bx) \geq r ) \leq \Pr ( E_n (\bx;r)^c )$, so that
\begin{equation}
\label{0830a}
 \Pr (R_n (\bx) > r) \leq \Pr \left( \bigcup_{i:1 \leq i \leq I,
\diam(C_i(\bx) \cap (0,1)^d) \geq r } \{ \U_n \cap C^+_i (\bx; r) = \emptyset \} \right) ,\end{equation}
with the convention that an empty union is empty.
Suppose $d \in (0,d^{1/2}]$.
For any $i$ with $\diam(C_i(\bx) \cap (0,1)^d) \geq r$, we can by convexity
choose a (non-random) $\bz \in C_i(\bx) \cap (0,1)^d$ at distance $r/2$ from $\bx$.
 Then (since $\frac{1}{4} < \frac{1}{2} \sin \frac{\pi}{12}$)
we have that $B (\bz ; r/4) \cap (0,1)^d$ is contained in
$C_i^+ (\bx;r)$ and,
since $r \leq d^{1/2}$,
 has $|B (\bz ; r/4) \cap (0,1)^d| \geq C r^d$ for some
$C \in (0,\infty)$ depending only on $d$.
Hence for  any $i$ with $\diam(C_i(\bx) \cap (0,1)^d) \geq r$,
\begin{equation}
\label{0830b} \Pr ( \U_n \cap C^+_i (\bx; r) = \emptyset )
\leq \Pr ( \U_n \cap B (\bz ; r/4) \cap (0,1)^d = \emptyset )
\leq (1-Cr^d)^n, \end{equation}
for some $C \in (0,\infty)$ depending only on $d$.
Applying Boole's inequality in (\ref{0830a}),
using (\ref{0830b}),
and noting that
  $1-x \leq \e^{-x}$ for any $x \geq 0$,
we have that
there are constants $C, C' \in (0,\infty)$,
  depending only on $d$,
such that
for all 
$r>0$ and $n \in \N$
\[ \sup_{\bx \in (0,1)^d} \Pr ( R_n(\bx) > r  ) \leq C' \exp (-C n r^d ).\]
Hence for $\beta>0$ and $n \in\N$, setting $s = C n r^{d/\beta}$, 
\begin{align*} \sup_{\bx \in (0,1)^d} \Exp [ R_n(\bx)^\beta ] & = 
 \sup_{\bx \in (0,1)^d} 
 \int_0^\infty \Pr ( R_n (\bx) > r^{1/\beta} ) \ud r
\leq C' \int_0^\infty \exp (-Cn r^{d/\beta} ) \ud r \\
& \leq 
C' n^{-\beta/d} \int_0^\infty s^{(\beta/d)-1} \exp (-s) \ud s
= C' n^{-\beta/d} \Gamma(\beta/d),
\end{align*}
using Euler's Gamma integral (see e.g.~6.1.1 in \cite{as})
for the last equality. $\qed$\\

Now we can complete the proof  of Lemma \ref{kkjj}.\\

\noindent
{\bf Proof of Lemma \ref{kkjj}.} 
With $R_n(\bx)$ as defined at (\ref{tttt}), we claim that 
\begin{equation}
\label{0830c}
\diam (V_n(\bx)) \leq  2 R_n(\bx)
\end{equation}
for  
all $\bx \in (0,1)^d$ and all $n \in \N$. 
Thus for $\beta>0$, $\Exp[ (\diam(V_n(\bx)))^\beta ] 
\leq C \Exp [ (R_n(\bx))^{\beta} ] \leq C' n^{-\beta/d}$,
by Lemma \ref{kkjjx}, proving Lemma \ref{kkjj}. 

To verify the claim (\ref{0830c}),
suppose that $\by \in (0,1)^d$ lies at distance
$s > r= R_n(\bx)$ from $\bx$. Then we can choose $i$
such that $\by \in C_i(\bx)$, so clearly
$\diam ( C_i(\bx) \cap (0,1)^d ) > r$
and 
$\diam ( C_i(\bx;r) \cap (0,1)^d ) =r$.
By definition of $R_n(\bx)$ we must have some point
of $\U_n \cap C_i^+(\bx;r)$; but then this
point lies closer to $\bx$ than $\by$ does,
so $\by$ is not in the Voronoi cell $V_n(\bx)$. Thus
$\sup_{\by \in V_n(\bx)} \| \bx -  \by \| \leq R_n (\bx)$.
Then the triangle inequality implies the result.
$\qed$\\

Next we establish the results that we will need later to control the changes in the ONG
on re-sampling a vertex.
Let $D \subset \R^d$ be a measurable, non-null convex region
and let $\bx \in D$.
Let $(\bX_1,\bX_2,\ldots)$ be a sequence of independent uniform random points on $D$.
We use the notation $\X_n = ( \bX_1,\ldots,\bX_n)$,
and for $\bx \in D$ set $\X_n^\bx := (\bx,\bX_1,\ldots,\bX_n)$.
For a finite sequence $\X$ of points in $\R^d$
and two points $\bx, \by$ of $\X$,
let $E(\bx,\by; \X)$
denote the event that $(\bx,\by)$ is an edge
in the ONG on $\X$.
Let $\O^{d,\alpha}_\bx (D;n)$ denote the total power-weighted length, 
with weight exponent $\alpha>0$, of edges incident to $\bx$
in the ONG on sequence $\X_n^\bx$, i.e.
\begin{align}
\label{defroo} \O^{d,\alpha}_\bx (D; n) := \sum_{i=1}^n \1_{E(\bX_i,\bx;  \X_i^\bx)}
\| \bX_i - \bx \|^\alpha .\end{align}
In the special case $D=(0,1)^d$, we
will write $\bU_i$ for $\bX_i$,
$\U_n$ for $\X_n$ and $\U_n^\bx$ for $\X_n^\bx$, and we
 abbreviate notation to
\[ \O^{d,\alpha}_\bx ( n) := \O^{d,\alpha}_\bx ((0,1)^d; n) .\]

\begin{lemma}
\label{xcxcxc}
Let $d \in \N$.
Suppose 
 $\alpha>0$. There exists $C  \in (0,\infty)$ such that
\begin{align}
\label{ook}
\sup_{n \in \N}
\sup_{\bx \in (0,1)^d}
\Exp [ \O^{d,\alpha}_\bx (n) ] \leq C.\end{align}
  Moreover there exists $C  \in (0,\infty)$ such that
for any $m$, $n$ with $0 \leq m < n$
\begin{align}
\label{ooky}
\sup_{\bx \in (0,1)^d}
\Exp [ \O^{d,\alpha}_\bx (n) - \O^{d,\alpha}_\bx (m)] 
\leq C (m+1)^{-\alpha/d}.\end{align}
\end{lemma}
\proof Fix $d \in \N$. 
For $i \in \N$ and $\bx \in (0,1)^d$, set
\[ W_i := \O^{d,1}_\bx (i) - \O^{d,1}_\bx (i-1)
= \1_{E(\bU_i,\bx;\U_i^\bx)} \| \bU_i -\bx \|,\]
with the convention $\O^{d,1}_\bx (0):=0$. Thus $W_i$ is the length of the edge
from $\bU_i$ to $\bx$ in the ONG on $\U_n^\bx$, if such an edge exists, or zero
otherwise. Then for $n \in \N$
\begin{align}
\label{eek}
 \O^{d,\alpha}_\bx (n) = \sum_{i=1}^n W_i^\alpha. 
 \end{align}
 Let $i \geq 2$.
 Given $\{\bU_1,\ldots,\bU_{i-1}\}$,  $W_i > 0$ only if $\bU_i$ falls inside the
 Voronoi cell of $\bx$ with respect to $\{\bx,\bU_1,\ldots,\bU_{i-1}\}$,
 that is $V_{i-1}(\bx)$ 
 as defined
at (\ref{voron}).
  In addition, given that $\bU_i \in V_{i-1}(\bx)$ (an
 event of probability $|V_{i-1}(\bx)|$), we have $W_i \leq {\rm diam} (V_{i-1}(\bx))$.
So for $i \geq 2$
\begin{align}
\label{bob11}
 \Exp [ W_i^\alpha \mid \bU_1,\ldots,\bU_{i-1} ] =
 \Exp [ W_i^\alpha \1_{\{ \bU_i \in V_{i-1} (\bx) \}} \mid \bU_1,\ldots,\bU_{i-1} ] \nonumber\\
 \leq 
|V_{i-1}(\bx)| (\diam(V_{i-1}(\bx)))^\alpha 
\leq (\diam( V_{i-1}(\bx) ))^{d+\alpha} 
.\end{align}
  Then
taking expectations in (\ref{bob11}) we obtain
\begin{align}
\label{bob1}
 \sup_{\bx \in (0,1)^d}
 \Exp [ W_i^\alpha ] \leq \sup_{\bx \in (0,1)^d}
 \Exp [ (\diam ( V_{i-1}(\bx)))^{d+\alpha} ]
 \leq C (i+1)^{-1-(\alpha/d)},\end{align}
for some $C \in (0,\infty)$ and all $i \in \N$,
by Lemma \ref{kkjj}.
Then we obtain (\ref{ook}) by taking expectations in (\ref{eek}) and
using (\ref{bob1}). 
Similarly we obtain (\ref{ooky}), this time using the fact that
for $1 \leq m < n$
\[ \Exp[ 
 \O^{d,\alpha}_\bx (n) - \O^{d,\alpha}_\bx (m)] = \sum_{i=m+1}^n \Exp[ W_i^\alpha]
 \leq C \sum_{i=m+1}^\infty (i+1)^{-1-(\alpha/d)}, \]
 by (\ref{bob1}).
 This completes the proof. 
$\qed$\\ 
 
In addition to $\O^{d,\alpha}_\bx (D;n)$, we consider the related quantity
\[ \hat \O^{d,\alpha}_\bx (D;n) := \sum_{i = 2  }^n
\1_{E(\bX_i,\bx; \X^\bx_i)}
(d (\bX_i; \{ \bX_1,\ldots,\bX_{i-1} \} ) )^\alpha;\]
that is, the total weight of the edges in the ONG on $\X_n$ from
those points that would be joined to $\bx$ in the ONG on $\X_n^\bx$.
In the case $D=(0,1)^d$, we use the abbreviation
\[ \hat \O^{d,\alpha}_\bx (n) := \hat \O^{d,\alpha}_\bx ((0,1)^d;n).\]
The following result is analogous to Lemma \ref{xcxcxc}.
\begin{lemma}
\label{lemhat}
Let $d \in \N$.
Suppose 
 $\alpha>0$. There exists $C 
 \in (0,\infty)$ such that
\begin{align}
\label{ookhat}
\sup_{n \in \N}
\sup_{\bx \in (0,1)^d}
\Exp [ \hat \O^{d,\alpha}_\bx (n) ] \leq C.\end{align}
  Moreover there exists $C  \in (0,\infty)$ such that
for any $m$, $n$ with $0 \leq m < n$
\begin{align}
\label{ookyhat}
\sup_{\bx \in (0,1)^d}
\Exp [ \hat \O^{d,\alpha}_\bx (n) - \hat \O^{d,\alpha}_\bx (m)] \leq C (m+1)^{-\alpha/d}.\end{align}
\end{lemma}
\proof
The proof is similar to that of
Lemma \ref{xcxcxc}. For $i \in\{2,3,4,\ldots\}$ set
\[ \hat W_i = \hat \O^{d,\alpha}_\bx (i)-\hat \O^{d,\alpha}_\bx (i-1)
= \1_{E( \bU_i, \bx;\U_i^\bx ) } d ( \bU_i ; \{ \bU_1,\ldots, \bU_{i-1} \} ),\]
where we take $\hat \O^{d,\alpha}_\bx (1):=0$.
Then $\hat W_i>0$ only if $\bU_i \in V_{i-1} (\bx)$.
Given that $\bU_i \in V_{i-1} (\bx)$, it follows from the triangle inequality that
\[ \hat W_i \leq d ( \bU_i ; \{ \bU_1,\ldots,\bU_{i-1} \} )
\leq d (\bU_i ; \{\bx\} ) + d (\bx ; \{ \bU_1,\ldots,\bU_{i-1} \} )
\leq C \diam ( V_{i-1} (\bx)) ,\]
for some $C \in (0,\infty)$
depending only on $d$.
It follows that there exists $C \in (0,\infty)$ such that
for all $i \in \N$
\[  \Exp [ \hat W_i^\alpha \mid \bU_1, \ldots,  \bU_{i-1} ] 
\leq C | V_{i-1} (\bx) | ( \diam ( V_{i-1} (\bx)) )^\alpha
\leq C ( \diam ( V_{i-1} (\bx)))^{d+\alpha} .\]
Thus by Lemma \ref{kkjj}, for some $C \in (0,\infty)$ and all $i \in \N$,
\[ \sup_{\bx \in (0,1)^d}  \Exp [ \hat W_i^\alpha ] \leq C
 (i+1)^{-1-(\alpha/d)},\]
and the lemma follows. $\qed$\\ 

The remaining results of this section will be used 
later to convert between Poisson and binomial
results.
The first 
is a technical lemma.

\begin{lemma}
\label{lem111}
Suppose $\beta \geq 1$ and $x>0$. Then,
\begin{align}
\label{111c} -\frac{1}{\beta} x^{1-\beta} \exp (-x^\beta ) \leq 
\int_0^x \exp(- t^\beta) \ud t - \Gamma ( 1+(1/\beta) ) \leq 0
.\end{align}
\end{lemma}
\proof
Suppose $\beta \geq 1$ and $x >0$. We have
\begin{equation}
\label{int0}
\int_0^x \exp(- t^\beta) \ud t
=
\int_0^\infty \exp(- t^\beta) \ud t - \int_x^\infty \exp(- t^\beta) \ud t .\end{equation}
We deal with each integral on the right-hand side of
(\ref{int0}) separately, using
the change of variable $y = t^\beta$.
By Euler's Gamma integral
(see e.g.~6.1.1 in \cite{as}) we have
\begin{equation}
\label{int1}
 \int_0^\infty \e^{-t^\beta} \ud t =
\frac{1}{\beta} \int_0^\infty y^{(1/\beta)-1} \e^{-y} \ud y
= \frac{1}{\beta} \Gamma(1/\beta) = \Gamma (1+(1/\beta)).
\end{equation}
For the second integral on the right-hand side
 of (\ref{int0}) we have
\begin{align}
\label{qwerty}
0 \leq \int_x^\infty \e^{- t^\beta} \ud t 
= \frac{1}{\beta} \int_{x^\beta}^\infty y^{(1/\beta)-1} \e^{-y} \ud y 
\leq \frac{1}{\beta} (x^\beta)^{(1/\beta)-1} \int_{x^\beta}^\infty \e^{-y} \ud y 
=  \frac{1}{\beta} x^{1-(1/\beta)}   \e^{-x^\beta} .\end{align}
Then from (\ref{int0}) with (\ref{int1}) and (\ref{qwerty})
we obtain (\ref{111c}).
$\qed$\\ 

To deduce the Poisson parts of 
Theorems \ref{ongvar} and  \ref{onngthm} 
we will need some 
estimates of
incremental
 expectations,
improving upon those in Section 3 of \cite{ong}. For $n \in \N$ set
\[ Z_n := \O^{d,1}(\U_n)-\O^{d,1}(\U_{n-1}),\]
taking $\O^{d,1}(\U_0) :=0$. Thus $Z_n$ is the gain in length on addition
of the $n$th point in the ONG on $(\bU_1,\bU_2,\ldots)$. Then for $n \in\N$
\begin{align}
\label{qqqz}
 \O^{d,\alpha} (\U_n) = \sum_{i=1}^n Z_i^\alpha.\end{align}
Note that (\ref{qqqz}) with (\ref{aaa}) below implies that for $\alpha \in (0,d)$
\begin{align}
\label{zzzf}
 \Exp [ \O^{d,\alpha} (\U_n) ] = \frac{d}{d-\alpha} v_d^{-\alpha/d} \Gamma(1+(\alpha/d))
n^{1-(\alpha/d)} + O ( \max \{ 1, n^{1-(\alpha/d)-(1/d)+\eps} \}),\end{align}
for any $\eps>0$, which improves upon the $o(n^{1-(\alpha/d)})$ error term
implicit in Theorem 2.1(i) of \cite{ong}. 

\begin{lemma}
\label{xp}
Suppose $d \in \N$ and $\alpha \in (0,d]$.
Then for $n \in \N$  
\begin{align}
\label{aaa}
 \Exp [ Z_n^\alpha ] =  v_d^{-\alpha/d} \Gamma(1+(\alpha/d)) n^{-\alpha/d}
+ h(n),\end{align}
where $h(n) = O(n^{-(\alpha/d)-(1/d)+\eps})$ as $n \to \infty$, for any $\eps>0$.
\end{lemma}
\proof 
Let $d\in\N$. For $r >0$ and $\bx \in (0,1)^d$, set $A(\bx;r)
:= |(0,1)^d \cap B (\bx;r)|$.
For $n\geq 2$,
\[ \Pr ( Z_n^\alpha > z \mid \bU_n ) = \Pr ( \{ \bU_1,\ldots,\bU_{n-1} \} \cap B( \bU_n ; z^{1/\alpha} )
= \emptyset \mid \bU_n)
 = ( 1 - A(\bU_n;z^{1/\alpha}) )^{n-1} .\]
For $r>d^{1/2}$, $A(\bx;r) = 1$ for all $\bx \in (0,1)^d$.
Then for $\bU_n \in (0,1)^d$,
\begin{align}
\label{hrnts}
 \Exp [ Z_n^\alpha \mid \bU_n ] = \int_0^\infty \Pr ( Z_n > z^{1/\alpha} \mid \bU_n ) \ud z
= \int_0^{d^{\alpha/2}} ( 1- A( \bU_n ; z^{1/\alpha} ) )^{n-1} \ud z .\end{align}
Fix $\eps \in (0,1/d)$ small. For all $n$ large enough so that
$n^{\eps-(1/d)} < 1/2$, 
let $S_n$ denote the region $[n^{\eps-(1/d)}, 1-n^{\eps-(1/d)}]^d$.
For 
$\bx=(x_1,x_2,\ldots,x_d) \in (0,1)^d$ let
$m(\bx):=\min \{x_1,\ldots,x_d,1-x_1,\ldots,1-x_d\}$, i.e.~the shortest distance from $\bx$ to the boundary
of $(0,1)^d$.
Consider $\bx \in S_n$. For $0 < r \leq m(\bx)$,  $A(\bx;r) = v_d r^d$, and for $r < d^{1/2}$, 
$C r^{d} \leq A(\bx;r) \leq v_d r^{d}$ for some $C \in (0,v_d)$ depending only on $d$. 
Thus from (\ref{hrnts})
\begin{align}
\label{111a}
 \Exp[ Z_n^\alpha \mid \bU_n \in S_n] 
 \geq  \int_0^{m(\bU_n)^\alpha} (1-v_d z^{d/\alpha})^{n-1} \ud z 
 \geq  \int_0^{n^{\eps\alpha-(\alpha/d)}} (1-v_d z^{d/\alpha})^{n-1} \ud z
,
\end{align}
since $m(\bU_n) \geq n^{\eps-(1/d)}$ for $\bU_n \in S_n$.
For $x>0$ Taylor's Theorem with Lagrange remainder implies that
$e^{-x} = 1 -x + C x^2$ where $C \in [0,1/2]$, so for $z < n^{\eps\alpha-(\alpha/d)}$ and $n$ large enough,
we have that
\begin{align*} \left( 1-v_d z^{d/\alpha} \right)^{n-1} & \geq 
\left( \exp ( -v_d z^{d/\alpha})
 -\frac{1}{2} v_d^2 z^{2d/\alpha} \right)^n \\
 & =  \exp ( -v_d n z^{d/\alpha} ) \left( 1- \frac{1}{2} v_d^2 z^{2d/\alpha} \exp ( v_d z^{d/\alpha})  
\right)^n \\
& \geq  \exp ( -v_d n z^{d/\alpha} ) \left( 1 + O \left( n^{2d\eps-2} \exp(v_d n^{d\eps-1})
\right)
\right)^n \\
& =  \exp ( -v_d n z^{d/\alpha} ) ( 1 + O(n^{2d\eps-1}) ),
 \end{align*}
 as $n \to \infty$, since $\eps < 1/d$. So from (\ref{111a}) we have that 
\begin{align}
\label{111b}
 \Exp[ Z_n^\alpha \mid \bU_n \in S_n] 
\geq (1+O(n^{2d\eps-1})) \int_0^{n^{\eps\alpha-(\alpha/d)}} \exp ( -v_d n z^{d/\alpha} ) \ud z
.\end{align}
Now,  setting $s = (v_d n)^{\alpha/d} z$, for $\alpha \in (0,d]$
\begin{align}
\label{qwert}
 \int_0^{n^{\eps\alpha-(\alpha/d)}} \exp ( -v_d n z^{d/\alpha} ) \ud z
& =  (nv_d)^{-\alpha/d} \int_0^{v_d^{\alpha/d} n^{\eps\alpha}} \exp ( - s^{d/\alpha} ) \ud s \nonumber\\
& = (n v_d)^{-\alpha/d} \Gamma (1+(\alpha/d)) + O ( \exp ( -v_d n^{\eps d} ) ),\end{align}
using (\ref{111c}) for the final equality.
So we obtain from (\ref{111b}) and (\ref{qwert})
that for $\eps>0$
\[ \Exp[ Z_n^\alpha \mid \bU_n \in S_n] \geq 
(n v_d)^{-\alpha/d}  \Gamma (1+(\alpha/d)) + O (n^{2d\eps-1-(\alpha/d)}).\] 
For the upper bound, using the fact that $1-x \leq \e^{-x}$ for $x \in (0,1)$
we have from (\ref{hrnts})
\begin{align}
\label{111d}
 \Exp[Z_n^\alpha \mid \bU_n \in S_n] =  
\int_0^{d^{\alpha/2}} (1-A(\bU_n;z^{1/\alpha}))^{n-1} \ud z \nonumber\\
 \leq  \int_0^{n^{\eps \alpha - (\alpha/d)}} \exp(-v_d (n-1) z^{d/\alpha}) \ud z
 + \int_{n^{\eps \alpha - (\alpha/d)}}^\infty \exp(- C (n-1) z^{d/\alpha}) \ud z.
\end{align}
For $\alpha \in (0,d]$, the second term on the right-hand side of (\ref{111d})
is $O( \exp ( -C n^{\eps d}))$, using (\ref{qwerty})
with $t= (C(n-1))^{\alpha/d} z$,  $\beta = d/\alpha$, and $x = (C(n-1))^{\alpha/d} n^{\eps \alpha - (\alpha/d)}$.
Also,
the first term on the right-hand side of (\ref{111d})
is bounded by
\[ 
\exp ( v_d n^{\eps d - 1}) \int_0^{n^{\eps \alpha - (\alpha/d)}} \exp(-v_d n z^{d/\alpha}) \ud z
= (nv_d)^{-\alpha/d} \Gamma (1+(\alpha/d)) + O( n^{\eps d - 1 -(\alpha/d)}),\]
by (\ref{qwert}).
So from (\ref{111d}), for the upper bound we obtain
\begin{align*} 
\Exp[Z_n^\alpha \mid \bU_n \in S_n]  \leq  
(n v_d)^{-\alpha/d} \Gamma (1+(\alpha/d))
+ O (n^{d\eps-1-(\alpha/d)}) .
\end{align*}
Combining the upper and lower bounds we have
\begin{align}
\label{1123s}
 \Exp[Z_n^\alpha \mid \bU_n \in S_n]  =  (n v_d)^{-\alpha/d} \Gamma (1+(\alpha/d))
+ O (n^{d\eps-1-(\alpha/d)}),\end{align}
for $\alpha \in (0,d]$ and $\eps$ small enough.
Now consider $\bx \in (0,1)^d \setminus S_n$. Here 
$Cr^d \leq A (\bx;r) \leq v_d r^d$ for $r<d^{1/2}$, and
by similar arguments to above,
we obtain
\begin{align}
\label{1123t}
 \Exp[Z_n^\alpha \mid \bU_n \notin S_n] = O(n^{-\alpha/d}).\end{align}
Since $\Pr(\bU_n \notin S_n)=O(n^{\eps-(1/d)})$, we obtain
from (\ref{1123s}) and (\ref{1123t}) that for any $\eps>0$
\begin{align*} \Exp [Z_n^\alpha] & = 
\Exp[Z_n^\alpha \mid \bU_n \in S_n] \Pr(\bU_n \in S_n)
+ \Exp[Z_n^\alpha \mid \bU_n \notin S_n] \Pr(\bU_n \notin S_n) \\
& = (n v_d)^{-\alpha/d} \Gamma (1+(\alpha/d))
+ O(n^{\eps-(\alpha/d)-(1/d)}),\end{align*}
and so
we have (\ref{aaa}). $\qed$

\section{Proof of Theorem \ref{ongvar}}
\label{vars}

The aim of this section is to prove the upper bounds on variances for $\O^{d,\alpha} (\U_n)$
and $\O^{d,\alpha} (\Po_\lambda)$ given in Theorem \ref{ongvar}.
The following martingale-difference
result is the key to the proof of the binomial parts of
Theorem \ref{ongvar}. 
Some extra work is then needed to derive the `Poissonized' version of the result.

\begin{lemma}
\label{mmmmm}
Let $d \in \N$ and $\alpha >0$. 
For each $n \in \N$, there exist mean-zero
random variables $D^{(n)}_i$, $i=1,2,\ldots,n$, such that:
\begin{itemize}
\item[(i)] $\sum_{i=1}^n D_i^{(n)} = \tO^{d,\alpha} (\U_n)$;
\item[(ii)] $\Exp[ D^{(n)}_i D^{(n)}_j ] =0$ for $i \neq j$;
\item[(iii)] there exists $C \in (0,\infty)$ such that
$\Exp[ (D^{(n)}_i)^2 ] \leq C i^{-2\alpha/d}$ for all $n$, $i$.
\end{itemize}
\end{lemma}

Before proving the lemma,
 we introduce some more notation.
For $n \in \N$, let $\F_n$ denote the $\sigma$-field
generated by $\U_n$.
Let $\F_0$ denote the trivial $\sigma$-field. 
For ease of notation during this proof, set 
$Y_n = \tO^{d,\alpha} (\U_n)$. 
Then we can write for $n\in\N$
\[ Y_n = \sum_{i=1}^n D_i^{(n)},\]
where for $i \in \{1,2,\ldots,n\}$
\begin{align}
\label{din}
 D^{(n)}_i = \Exp[ Y_n \mid \F_i ] - \Exp[ Y_n \mid \F_{i-1}],\end{align}
 and for fixed $n$
 the $D_i^{(n)}$, $i=1,\ldots,n$ are martingale differences, and hence orthogonal (see e.g.~Chapter 12 of \cite{will}). This establishes
 parts (i) and (ii) of the lemma. It remains  to
estimate $\Exp[(D_i^{(n)})^2]$.
Given $\U_n = (\bU_1,\ldots,\bU_n)$, for $i \in \{1,\ldots,n\}$
let $\bU_i'$ be an independent copy of
$\bU_i$ (independent of $\bU_1,\bU_2,\ldots$) and
set 
\[ \U_n^i := (\bU_1,\ldots,\bU_{i-1}, \bU_i', \bU_{i+1},\ldots, \bU_n ),\]
so $\U_n^i$ is $\U_n$  with the $i$th member of the sequence independently
re-sampled. Define
\[ \Delta^{(n)}_i := \tO^{d,\alpha} ( \U^i_n ) - \tO^{d,\alpha}
(\U_n ) = \O^{d,\alpha} ( \U^i_n ) - \O^{d,\alpha}
(\U_n ),\]
the change in $Y_n$ on re-sampling the point $\bU_i$. 
Then it is the case that
\[ D^{(n)}_i = - \Exp[ \Delta^{(n)}_i \mid \F_i ] .\]
We split $\Delta^{(n)}_i$ into six
 components as follows. Let $\Delta^{(n)}_{i,1}$ be the
weight of the  edge from $\bU_i$ in the ONG on $\U_n$, and let $\Delta^{(n)}_{i,2}$ be the
weight of the   edge from $\bU'_i$ in the ONG on $\U^i_n$. Let $\Delta^{(n)}_{i,3}$ be the
total weight of the edges incident to $\bU_i$ in the ONG on $\U_n$, and let $\Delta^{(n)}_{i,4}$ be the
total weight of the edges incident to $\bU'_i$ in the ONG on $\U^i_n$. 
Let $\Delta^{(n)}_{i,5}$ be the total weight of edges in the
ONG on $(\bU_1,\ldots,\bU_{i-1},\bU_{i+1},\ldots,\bU_n)$ from points
in $\U_n$ that are joined to $\bU'_i$ in the ONG on $\U^i_n$.
Let $\Delta^{(n)}_{i,6}$ be the total weight of edges in the
ONG on $(\bU_1,\ldots,\bU_{i-1},\bU_{i+1},\ldots,\bU_n)$ from points
in $\U_n$ that are joined to $\bU_i$ in the ONG on $\U_n$.
Then
\[ \Delta^{(n)}_i = \Delta^{(n)}_{i,2} + \Delta^{(n)}_{i,4} +\Delta^{(n)}_{i,6}
-\Delta^{(n)}_{i,1}- \Delta^{(n)}_{i,3} - \Delta^{(n)}_{i,5}.\]
The next result will be crucial for the proof of Lemma \ref{mmmmm}.

\begin{lemma}
\label{clem}
For any $\alpha>0$
there exists $C \in (0,\infty)$ such that
for all $\ell \in\{1,\ldots,6\}$
\begin{align}
\label{claim}
\Exp [ (\Exp [ \Delta_{i,\ell}^{(n)} \mid \F_i ])^2 ] \leq C i^{-2\alpha/d} ,\end{align}
for all $n\in \N$ and $i\in \{1,\ldots,n\}$. \end{lemma}
\proof
First consider $\ell \in \{1,2\}$.
By the conditional Jensen's inequality,
\[ \Exp [ (\Exp [ \Delta^{(n)}_{i,\ell} \mid \F_i ])^2 ] 
\leq \Exp [ \Exp [ (\Delta^{(n)}_{i,\ell})^2 \mid \F_i ] ]
= \Exp [ (\Delta^{(n)}_{i,\ell})^2 ].\]
For $\ell \in \{1,2\}$, we have from Lemma 3.1 in \cite{ong} (cf (\ref{aaa}) above) that
for $\alpha>0$,  
$\Exp [ (\Delta^{(n)}_{i,\ell})^2 ] = \Exp [ Z_i ^{2\alpha} ] \leq C i^{-2\alpha/d}$
for all $i, n$.
 Thus
 for $\ell \in \{1,2\}$ there
is a constant $C \in (0,\infty)$ such that (\ref{claim}) holds
for all $i$ and $n$. 

Now consider $\ell\in \{3,4\}$.
For $i \in \N$, let $V_i:=V_{i-1}(\bU_i)$ be the Voronoi cell of $\bU_i$ with respect to 
$\{\bU_1,\ldots,\bU_{i-1},\bU_i\}$.
Similarly, let $V'_i:=V_{i-1}(\bU'_i)$ be the Voronoi cell of $\bU'_i$ with respect to
$\{\bU_1,\ldots,\bU_{i-1},\bU'_i\}$. 

By convexity, there exists a $d$-cube of side length
at most $2 \diam (V_i)$ which contains $V_i$ and
also lies inside $(0,1)^d$. 
Let $B_i$ denote a
minimal-volume such cube.

Points of $\{\bU_{i+1}, \ldots, \bU_n\}$ that 
fall outside of $V_i$ can never be joined to $\bU_i$ and can only serve to 
decrease the total weight incident to $\bU_i$
(by shrinking the
subsequent Voronoi cells). Hence
removing any point
of $\{\bU_{i+1},\ldots,\bU_n\}$ that falls
 outside $V_i$ (and in particular any that falls
 outside $B_i$)
 can only increase
the total weight of edges incident to $\bU_i$. Moreover,
$\{ \bU_1,\ldots, \bU_{i-1}\}$ necessarily
lie outside $V_i$ and their removal can only increase the total
weight incident to $\bU_i$. In other words,
for any $j \geq i+1$ and any
subsequence $\U'_j$ of $\U_j$ containing $\bU_i$ and $\bU_j$,
we have $E(\bU_j,\bU_i;\U_j) \subseteq E(\bU_j , \bU_i ; \U'_j )$, and
$\Pr ( E (\bU_j, \bU_i ; \U_j )) = 0$ for any $\bU_j \notin V_i$
and in particular any $\bU_j \notin B_i$.

It follows that
\[ \Delta_{i,3}^{(n)} = \sum_{j=i+1}^n \1_{E(\bU_j,\bU_i; \U_j)} \| \bU_j - \bU_i \|^\alpha
\leq  \sum_{j: i+1 \leq j \leq n,\bU_j \in B_i} \1_{E(\bU_j,\bU_i; \U_{j,i}  )} \| \bU_j - \bU_i \|^\alpha ,\]
where $\U_{j,i}$ is the subsequence of $(\bU_i,\ldots,\bU_j)$ consisting only of those
points in $B_i$.
So in particular, given $\F_i$, $\Delta_{i,3}^{(n)}$
is stochastically dominated by $\O_{\bU_i}^{d,\alpha} ( B_i ; N )$ where
$N \sim {\rm Bin} (n-i,|B_i|)$ is the number of points
of $\{\bU_{i+1},\ldots,\bU_n\}$ that fall  in $B_i$. 
(Recall the definition of $\O^{d,\alpha}_\bx (D; n)$
from (\ref{defroo}).) 
We thus have that, given $\F_i$,
$\Delta^{(n)}_{i,3}$
is stochastically dominated by 
\[ \O^{d,\alpha}_{\bU_i} (B_i; n) 
\eqd |B_i|^{\alpha/d} \O^{d,\alpha}_\bx (  n),\]
by scaling, for some
 $\bx \in (0,1)^d$.
 Since $|B_i|\leq C (\diam (V_i))^d$, we have
 in particular that for all $n \in\N$ and $i \in \{ 1,\ldots,n \}$
\[ \Exp [ \Delta^{(n)}_{i,3} \mid \F_i ] \leq  C (\diam (V_i))^{\alpha}
\sup_{\bx \in (0,1)^d} \Exp [ \O^{d,\alpha}_\bx (  n) ]
\leq  C (\diam(V_i))^{\alpha} ,\]
by (\ref{ook}). Thus by Lemma \ref{kkjj},
for all $n \in \N$ and $i\in \{1,\ldots,n\}$,
\begin{align}
\label{bbbb1}
 \Exp[ (\Exp [ \Delta^{(n)}_{i,3} \mid \F_i ])^2 ] \leq 
C \Exp [   (\diam(V_i))^{2\alpha} ] \leq C i^{-2\alpha/d}.\end{align}
Similarly, $\Exp [ \Delta^{(n)}_{i,4} \mid \F_i ] \leq C \Exp [ (\diam(V'_i))^{\alpha} \mid \F_i ]$
so that, by the conditional Jensen's inequality,
\begin{align}
\label{bbbb2}
 \Exp [ (\Exp [ \Delta^{(n)}_{i,4} \mid \F_i ])^2 ] \leq 
  C \Exp [ (\diam(V'_i))^{2\alpha}]
\leq C i^{-2\alpha/d},\end{align}
for all $n \in \N$ and $i\in\{1,\ldots,n\}$ by Lemma \ref{kkjj} once more, since
$V'_i \eqd V_i$. Thus
from (\ref{bbbb1}) and (\ref{bbbb2}) we verify
the $\ell \in \{3,4\}$ cases of (\ref{claim}).

Finally consider $\ell \in \{5,6\}$. 
Recall that $\U_{j,i}$ is the subsequence
of $(\bU_i,\ldots,\bU_j)$ consisting only of those points in $B_i$.
By the argument above for $\Delta_{i,3}^{(n)}$, we
have that
 \begin{align}
\label{0831a}
 \Delta_{i,6}^{(n)} 
 = \sum_{j=i+1}^n \1_{E(\bU_j,\bU_i; \U_j)} 
\left( d(\bU_j; \{ \bU_1,\ldots,\bU_{j-1}\} \setminus
\{ \bU_i \} ) \right)^\alpha \nonumber\\
 \leq  \sum_{j: i+1 \leq j \leq n, \bU_j \in B_i} 
 \1_{E(\bU_j,\bU_i; \U_{j,i})} 
 \left(d(\bU_j; \{ \bU_1,\ldots,\bU_{j-1}\} \setminus
\{ \bU_i \} ) \right)^\alpha.\end{align}
List the points of $\U_{j,i}$  
 in order of increasing mark (index)
as $(\bU_i,\bU_{j_1},\ldots,\bU_{j_s})$. 
For $j \geq j_1 +1$, observe that 
removing points outside $B_i$ can only increase the distance
from $\bU_j$ to its nearest
neighbour amongst $\{ \bU_1,\ldots, \bU_{j-1} \} \setminus \{ \bU_i\}$,
since we know $\bU_{j_1} \in  B_i$.
Thus we have that for $j  \geq j_1 +1$
\begin{equation}
\label{0831b}
 d ( \bU_{j} ; \{ \bU_1,\ldots,\bU_{j-1} \} \setminus \{ \bU_i\} )
\leq 
d ( \bU_{j} ; ( \{ \bU_1,\ldots,\bU_{j-1} \} \setminus \{ \bU_i\} ) \cap B_i ).\end{equation}
Then from (\ref{0831a}) and (\ref{0831b}) we obtain
\begin{align}
\label{0831d}
 \Delta_{i,6}^{(n)}  & \leq \left( d( \bU_{j_1} ; \{ \bU_1,\ldots,\bU_{j_1-1} \} \setminus \{ \bU_i\}   ) \right)^\alpha
\nonumber\\
& ~~+
 \sum_{j: i+1 \leq j \leq n, \bU_j \in B_i} 
 \1_{E(\bU_j,\bU_i; \U_{j,i})} 
 \left(d(\bU_j; ( \{ \bU_1,\ldots,\bU_{j-1}\} \setminus
\{ \bU_i \} ) \cap B_i ) \right)^\alpha.\end{align}
To bound the length of the edge from $\bU_{j_1}$, we note that
any point $\by \in V_n(\bx)$ has $d ( \by; \{ \bU_1,\ldots, \bU_{n} \} ) \leq
2 \diam (V_n(\bx))$. Hence
\begin{equation}
\label{0831c}
 d ( \bU_{j_1} ; \{ \bU_1,\ldots,\bU_{j_1-1} \} \setminus \{ \bU_i\} )
\leq C \diam (V_i).\end{equation}
Recall the definition of $\hat \O^{d,\alpha}_\bx (D;n)$
from just above Lemma \ref{lemhat}.
Then from (\ref{0831d}) with   (\ref{0831c}),
 we have 
 that, given $\F_i$,
  $\Delta_{i,6}^{(n)}$ is stochastically
dominated by
\[ \hat \O^{d,\alpha}_{\bU_i} ( B_i ; n)  + C ( \diam (V_i))^\alpha
\eqd |B_i|^{\alpha/d} \hat \O^{d,\alpha}_\bx ( n) + C ( \diam (V_i))^\alpha,\]
for some $\bx \in (0,1)^d$. Taking expectations, we obtain from Lemma
\ref{lemhat} that
\[ \Exp [ \Delta_{i,6}^{(n)} \mid \F_i ] \leq C |B_i|^{\alpha/d} + C ( \diam (V_i))^\alpha
\leq C' ( \diam (V_i))^\alpha .\]
Then by Lemmas
  \ref{kkjj} we obtain
\[ \Exp[ (\Exp [ \Delta_{i,6}^{(n)} \mid \F_i ])^2 ]   
\leq C \Exp [ (\diam (V_i))^{2\alpha} ] \leq C' i^{-2\alpha/d},\]
for all $n, i$. A similar argument holds for $\Delta_{i,5}^{(n)}$,
and thus verifies the $\ell \in \{5,6\}$ cases of (\ref{claim}).
This completes the proof of the lemma.
 $\qed$\\ 
 
\noindent
{\bf Proof of Lemma \ref{mmmmm}.}
With $D_i^{(n)}$ given by (\ref{din}),
parts (i) and (ii) of the lemma are immediate, 
as described above.
The Cauchy--Schwarz inequality and (\ref{claim}) imply 
 \[ \Exp [( D^{(n)}_i)^2] = \Exp [ (\Exp[ \Delta^{(n)}_i \mid \F_i ])^2 ]
 = \Exp \left[ \left( \sum_{\ell=1}^6 (-1)^\ell \Exp [ \Delta^{(n)}_{i,\ell} \mid \F_i ] \right) ^2\right]
  \leq C i^{-2\alpha/d},\]
 for all $n, i$. This yields part (iii) of the lemma. $\qed$\\ 
  
To deduce the Poisson version of Theorem \ref{ongvar}, and later Theorem \ref{onngthm},
 we prove the following 
series of
lemmas. 
\begin{lemma}
\label{poisss}
Let $N(n)$ be a Poisson random variable with mean $n \geq 1$. For $\beta \in [ 0,1)$,
\begin{align}
\label{vv1}
 \Var [ N(n)^{1-\beta} ] \leq C n^{1-2\beta}; \\
 \label{vv2}
 \Exp [ ( N(n)^{1-\beta} - n^{1-\beta} )^2 ] \leq C n^{1-2\beta}; \\
 {\rm and }~~~
 \label{vv3}
 \Exp [ (\log (1+N(n)) - \log (1+n))^2 ] \leq C n^{-1} ;\end{align}
 for some $C \in (0,\infty)$ and all $n \geq 1$. 
\end{lemma}
\proof Let $n \geq 1$.
First we prove (\ref{vv1}), (\ref{vv2}). 
Let $\beta \in [0,1)$.
Set $K_n := N(n) -n$. Then
\begin{align}
\label{nev}
  N(n)^{1-\beta} = n^{1-\beta} (1 + n^{-1} K_n )^{1-\beta},\end{align}
where by
the Intermediate Value Theorem we have that
$(1 + n^{-1} K_n )^{1-\beta} = 1 + (1-\beta) (1+ H_n)^{-\beta} n^{-1} K_n$ for some
$H_n$ with $|H_n| \leq n^{-1}|K_n|$. Hence
\begin{align}
\label{nev1}
  N(n)^{1-\beta} - n^{1-\beta}  = n^{-\beta} (1-\beta) (1+H_n)^{-\beta} K_n,\end{align}
so that for $C \in (0,\infty)$
\begin{align}
\label{nev2}
\Exp [ ( N(n)^{1-\beta} - n^{1-\beta} )^2 ] 
= C n^{-2\beta} \Exp [ (1+H_n)^{-2\beta} K_n^2 ] .\end{align}
Let $A_n$ denote the event $\{ |K_n| < n^{3/4}\}$.
Then
\[ \Exp [ (1+H_n)^{-2\beta} K_n^2 ] = \Exp[ (1+H_n)^{-2\beta} K_n^2 \1_{A_n}]
+ \Exp[ (1+H_n)^{-2\beta} K_n^2 \1_{A_n^c}].\]
Here, by Cauchy--Schwarz,
\[ \Exp[ (1+H_n)^{-2\beta} K_n^2 \1_{A_n^c}] \leq (\Exp[ (1+H_n)^{-4\beta} K_n^4 ])^{1/2}
(\Pr (A_n^c))^{1/2} .\]
But by (\ref{nev1}), for $C \in (0,\infty)$,
$\Exp[ (1+H_n)^{-4\beta} K_n^4 ] =  C n^{4\beta} \Exp [  | N(n)^{1-\beta} - n^{1-\beta} |^4 ]$,
so that
\begin{align}
\label{nev3}
 \Exp[ (1+H_n)^{-2\beta} K_n^2 \1_{A_n^c}] \leq 
C n^{2\beta} ( \Exp [  | N(n)^{1-\beta} - n^{1-\beta} |^4 ]
 )^{1/2} (\Pr (A_n^c))^{1/2} ,\end{align}
which tends to zero as $n \to \infty$, by 
standard Chernoff-type
 Poisson tail bounds (see e.g.~Lemma 1.2 in \cite{penbook}). 
Also, given $A_n$, $|H_n| \leq n^{-1}|K_n| \leq n^{-1/4}$,
so that
\begin{align}
\label{nev4} \Exp[ (1+H_n)^{-2\beta} K_n^2 \1_{A_n}]
\leq C \Exp [ K_n^2 \mid A_n ] 
\leq C \Exp[ K_n^2] \Pr (A_n)^{-1} = C n \Pr(A_n)^{-1} \sim C n,\end{align}
as $n \to \infty$, by standard Poisson tail bounds. So from
(\ref{nev2}), (\ref{nev3}) and (\ref{nev4}) we obtain (\ref{vv2}).

Now from (\ref{nev1}) we have
\[ \Var [ N(n)^{1-\beta} ] = C n^{-2\beta}
\Var [ (1+H_n)^{-\beta} K_n ] \leq  C n^{-2\beta} \Exp [ (1+H_n)^{-2\beta} K_n^2 ]
.\]
Then from (\ref{nev3}) and (\ref{nev4}) we obtain (\ref{vv1}).

Finally, the Intermediate Value Theorem 
 implies that
\[ \log (1+N(n)) - \log (1+n) =  \log (1 + (1+n)^{-1} K_n ) = (1+n)^{-1} (1+ H_n)^{-1} K_n,\]
where, as before, $|H_n| \leq n^{-1} |K_n|$. Hence
\[ \Exp [ ( \log (1+N(n)) - \log (1+n))^2 ]
= (1+n)^{-2} \Exp [ (1+H_n)^{-2} K_n^2 ]. \]
Now  (\ref{nev4}) still holds with $\beta=1$, while instead of (\ref{nev3})
in this case we have 
\[ \Exp [ (1+H_n)^{-2} K_n^2 \1_{A_n^c} ]
\leq (1+n)^2 (\Exp [ (\log (1+N(n)) - \log (1+n))^4 ])^{1/2} (\Pr (A_n^c))^{1/2},\]
by Cauchy--Schwarz,
which again tends to zero as $n\to\infty$.
Thus we obtain (\ref{vv3}).
$\qed$

\begin{lemma}
Let $d \in \N$ and $\alpha >0$. Let $N(n)$ be a Poisson random variable with mean $n \geq 1$. Then
there exists $C \in (0,\infty)$ such that for all $n \geq 1$
\begin{align}
\label{pyy}
\Exp [ \Var [\O^{d,\alpha}(\U_{N(n)}) \mid N(n)]]
\leq C +  \sup_{1 \leq m \leq 2n} \Var [ \O^{d,\alpha} (\U_{m}) ].\end{align}
\end{lemma}
\proof 
We have that
\[ \Exp [ \Var [\O^{d,\alpha} (\U_{N(n)}) \mid N(n) ]]
\leq \sup_{1 \leq m \leq 2n} \Var [ \O^{d,\alpha} (\U_{m} )]
+ C \Exp [ (N(n))^2 \1_{ \{ N(n) > 2n \}} ],\]
using the trivial bound that $\O^{d,\alpha} (\U_{N(n)}) \leq C N(n)$. By Cauchy--Schwarz, the last
term in the above display is bounded by a constant times 
\[ (\Exp [ (N(n))^4 ])^{1/2} (\Pr ( N(n) >2n))^{1/2},\]
which tends to $0$ as $n\to\infty$ by standard Poisson tail bounds.
So we obtain 
(\ref{pyy}).  
$\qed$

\begin{lemma} 
\label{nevv}
Let $d \in \N$ and $\alpha \in (0,d]$. 
For $n \in \N$, let $\mu_n := \Exp[ \O^{d,\alpha} (\U_n)]$. Let $N(n)$ be a Poisson
random variable with mean $n$. There exists $C \in (0,\infty)$ such that
 for all $n \geq 1$
 \[ \Exp [  ( \mu_{N(n)} - \mu_{\lfloor n \rfloor} )^2 ]
 \leq Cn^{1-(2\alpha/d)} .\]
 \end{lemma}
 \proof 
 Taking expectations in (\ref{qqqz}), 
 we have that for $n \in\N$
\[ \mu_n = \sum_{i=1}^n \Exp [Z_i^\alpha].\]
First suppose that $\alpha \in (0,d)$. By (\ref{aaa}) we have that, for integers $\ell$, $m$ with
$1 \leq \ell < m$,
\begin{align}
\label{ggjj}
 \mu_m-\mu_{\ell} =\sum_{i=\ell+1}^m \Exp[Z_i^\alpha] = \frac{d}{d-\alpha} v_d^{-\alpha/d} \Gamma(1+(\alpha/d))
(m^{1-(\alpha/d)} - \ell^{1-(\alpha/d)} ) \nonumber\\
+ \sum_{i=\ell+1}^m h(i) + O(m^{-\alpha/d}) + O(\ell^{-\alpha/d}).\end{align}
In particular, for $n \geq 1$
\begin{align}
\label{bow1}
 | \mu_{N(n)} -\mu_{ \lfloor n \rfloor} |  = C | N(n)^{1-(\alpha/d)} - n^{1-(\alpha/d)} | + \delta(n),\end{align}
where from (\ref{ggjj})
the random variable $\delta(n)$ satisfies
\begin{align}
\label{kk66}
 | \delta(n) | \leq 
    \sum_{i=\min \{ N(n),  \lfloor n \rfloor  \} }^{\max \{  N(n), \lfloor n \rfloor \}} | h(i) |
  + O ( \min \{n , 1+N(n) \}^{-\alpha/d} ) . \end{align}

On the other hand, for $\alpha = d$, this time
(\ref{aaa}) implies that for
$1 \leq \ell < m$,
\begin{align}
\label{ppiiuu}
 \mu_m-\mu_{\ell} =  v_d^{-1} 
(\log (1+m) - \log (1+\ell) )
+ \sum_{i=\ell+1}^m h(i) + O(m^{-1}) + O(\ell^{-1}).\end{align}
In particular, for $n \geq 1$, (\ref{ppiiuu}) gives
\begin{align}
\label{bow11}
 | \mu_{N(n)} -\mu_{ \lfloor n \rfloor} |  = C | \log (1+N(n)) - \log (1+n) |
 + \delta(n),\end{align}
 where again $\delta (n)$ satisfies (\ref{kk66}), now with $\alpha=d$.
 
 We now claim that for all $\alpha \in (0,d]$,
 $\delta(n)$ as defined by (\ref{bow1}) or (\ref{bow11})
 satisfies
 \begin{align}
\label{bbnnmm}
\Exp[ \delta(n)^2 ] = o (n^{1-(2\alpha/d)}) ,
~{\rm as}~ n\to\infty.\end{align}
Then in the case $\alpha \in (0,d)$,
(\ref{bbnnmm}) with (\ref{bow1}), (\ref{vv2})
and  Cauchy--Schwarz yields the lemma.
In the case $\alpha=d$, the result follows from
(\ref{bbnnmm}) with (\ref{bow11}), (\ref{vv3})
and Cauchy--Schwarz again.

It remains to prove the claim (\ref{bbnnmm}). We start 
from the fact that for $\alpha \in (0,d]$, $\delta(n)$
satisfies (\ref{kk66}).
Note that there exists $C \in (0,\infty)$
such that for $n \geq 1$
\begin{align}
\label{lapel} \Exp [ \min \{ n , 1+ N(n)\}^{-2\alpha/d } ]
\leq n^{-2\alpha/d} + \Exp [ (1+N(n))^{-2\alpha/d} ]
\leq C n^{-2\alpha/d} ,\end{align}
as can be proved
by standard Poisson tail estimates
as used elsewhere in the present paper
(cf Lemma \ref{binbin} for an analogous binomial result).

Now we deal with the main term
in (\ref{kk66}). 
We have that  for $n \geq 1$,
 $\alpha \in (0,d]$
and   any $\eta \in (0,1/2)$,
\begin{align}
\label{njkml}
 \sup_{m \in \N: | m - \lfloor n \rfloor|  \leq n^{(1/2)+\eta} }  
 \sum_{i=\min \{m,\lfloor n \rfloor\}}^{\max\{m,\lfloor n \rfloor\}}
  |h(i) |
\leq  ( 2 n^{(1/2)+\eta} +1)
 \sup_{ m \in \N: | m - \lfloor n \rfloor|  \leq n^{(1/2)+\eta} } | h(m) | ;
\end{align}
it follows from (\ref{njkml}) and
  Lemma \ref{xp} that for any $\eta \in (0,1/2)$,
$\eps>0$, there exists $C \in (0,\infty)$ such that for
all $n \geq 1$
\begin{align}
\label{kkkmmm}
\sup_{m\in \N: | m - \lfloor n \rfloor|  \leq n^{(1/2)+\eta} }  
\sum_{i=\min \{m,\lfloor n \rfloor\}}^{\max\{m,\lfloor n \rfloor\}}
  |h(i) |
\leq C n^{(1/2)+\eta-(1/d)-(\alpha/d)+\eps} .
\end{align}
In particular, this is $o(n^{(1/2)-(\alpha/d)})$ for  sufficiently small
$\eps$, $\eta$. 
Now we have, with $\eta>0$ as above, for $n \geq 1$
\begin{align}
\label{kkkmmml}
 \Exp \left[ \left(\sum_{i=\min\{N(n) , \lfloor n \rfloor\}}^{\max\{N(n) , \lfloor n \rfloor\}} 
 | h(i) | \right)^2 \right] \leq 
  \left( \sup_{m \in \N: |m-\lfloor n \rfloor| \leq n^{(1/2)+\eta} }  \sum_{i=\min \{m,\lfloor n \rfloor\}}
^{\max\{m,\lfloor n \rfloor\}} |h(i)| \right)^2 \nonumber\\
+ C  \Exp [ |N(n)-\lfloor n \rfloor|^2 \1_{\{|N(n)-\lfloor n \rfloor| > n^{(1/2)+\eta}\}} ],\end{align}
and by Cauchy--Schwarz
\[ \Exp [ |N(n)-\lfloor n \rfloor|^2 \1_{\{|N(n)-\lfloor n \rfloor| > n^{(1/2)+\eta}\}} ]
\leq (\Exp [ |N(n)-\lfloor n \rfloor|^4 ]  \Pr (|N(n)-\lfloor n \rfloor| > n^{(1/2)+\eta }))^{1/2}  ,\]
which is $o(n^{1-(2\alpha/d)})$
as $n \to \infty$, by standard Poisson tail bounds. Thus from (\ref{kkkmmml}),
(\ref{kkkmmm}), (\ref{lapel}), and Cauchy--Schwarz,
we 
verify (\ref{bbnnmm}).
$\qed$\\ 

\noindent
{\bf Proof of Theorem \ref{ongvar}.} First we prove the binomial parts 
of (\ref{a1}) and (\ref{a3}). By part (i) of 
Lemma \ref{mmmmm}, we have that
$\tO^{d,\alpha} (\U_n) = \sum_{i=1}^n D_i^{(n)}$ for each $n \in \N$.
 By the orthogonality of the $D_i^{(n)}$ (part (ii) of Lemma \ref{mmmmm})
  we have that for $n \in \N$
\[ \Var [\tO^{d,\alpha}(\U_n)] = \sum_{i=1}^n \Exp[(D_i^{(n)})^2] ,\]
which by part (iii) of Lemma \ref{mmmmm} yields the upper bounds as claimed. 

We now deduce the Poisson 
parts of (\ref{a1}) and (\ref{a3}).
For ease of notation, let 
$X_n := \O^{d,\alpha} (\U_n)$ and $\mu_n := \Exp[X_n]$. Then if $N(n)$ is Poisson with mean
$n \geq 1$, 
$\O^{d,\alpha} (\Po_n)$ has the distribution
of $X_{N(n)}$ and its expectation is $\Exp[ \mu_{N(n)} ]=:a_n$. 
Write
\begin{align}
\label{iiii}
\tO^{d,\alpha} (\Po_n) =
X_{N(n)} - a_n = (X_{N(n)} - \mu_{N(n)} ) + (\mu_{\lfloor n \rfloor} - a_n) + (\mu_{N(n)} -\mu_{\lfloor n \rfloor} ).\end{align}
Then $\Var [ \tO^{d,\alpha} (\Po_n)] = \Var [ (X_{N(n)} - \mu_{N(n)} ) + (\mu_{N(n)} -\mu_{\lfloor n \rfloor} ) ]$.
We have
\[ \Var [ X_{N(n)} - \mu_{N(n)} ] = \Exp [ \Var [ X_{N(n)} - \mu_{N(n)} \mid N(n) ]]
= \Exp [ \Var [ \O^{d,\alpha} ( \U_{N(n)}) \mid N(n)]].\]
By (\ref{pyy}) this is bounded by a constant times $\sup_{m \leq 2n} \Var [ \O^{d,\alpha} (\U_m)]$,
which,
using the binomial parts
 of (\ref{a1}) and (\ref{a3}),
is bounded by a constant times $n^{1-(2\alpha/d)}$ for  $\alpha \in (0,d/2)$ and
by a constant times $\log (1+n)$ for $\alpha =d/2$.
So
we have 
for $C \in (0,\infty)$ and $n \geq 1$
\begin{align}
\label{bow3}
 \Var [ X_{N(n)} - \mu_{N(n)} ] 
\leq 
\left\{
\begin{array}{ll} 
C n^{1-(2\alpha/d)} & \textrm{ if } \alpha \in (0,d/2); \\
 C \log (1+n) & \textrm{ if } \alpha=d/2 .\end{array}\right. \end{align}
The final term on the right-hand side
of (\ref{iiii}) satisfies Lemma \ref{nevv}.
So by (\ref{iiii}) with 
Lemma \ref{nevv},  (\ref{bow3}), 
and Cauchy--Schwarz,
 we obtain the Poisson parts of (\ref{a1}) and (\ref{a3}).
 $\qed$

\section{Proof of Theorem \ref{onngthm}}
\label{prf}

By Lemma \ref{mmmmm} we have that for $\alpha > d/2$, for all $n \in \N$
\[ \Var [ \tO^{d,\alpha} (\U_n) ] = \sum_{i=1}^n \Exp[(D_i^{(n)})^2] \leq C \sum_{i=1}^n i^{-2\alpha/d}
\leq C' < \infty.\]
 In order to show that $\tO^{d,\alpha} (\U_n)$ in fact
 converges, we employ a refinement of the
 martingale difference technique of Section \ref{vars}. First we need two more lemmas.
\begin{lemma}
\label{binbin}
Suppose $X \sim {\rm Bin} (n,p)$ for $n \in \N$ and
 $p \in (0,1)$. Then for any $\beta >0$
there exists $C \in (0,\infty)$ such that for all $n \in \N$
and all $p \in (0,1)$
\[ \Exp [ (1+X)^{-\beta}] \leq C (np)^{-\beta}.\]
\end{lemma}
\proof
We have that
\begin{align*} \Exp [ (1+X)^{-\beta} ] & \leq 
 ( 1 + (np/2))^{-\beta} + \Exp [ (1+X)^{-\beta} \1_{ \{ X < np/2 \}} ]
\\
& \leq 
C (n p)^{-\beta} + (\Exp [ (1+X)^{-2\beta} ]) ^{1/2} (\Pr ( X < np/2 ))^{1/2},\end{align*}
for some $C \in (0,\infty)$ and all $n \in \N$, $p \in (0,1)$, 
using Cauchy--Schwarz. But for $\beta
>0$, $(1+X)^{-2\beta} \leq 1$ a.s., so $\Exp [ (1+X)^{-2\beta} ] \leq 1$. Also, by standard
binomial tail bounds (see e.g. Lemma 1.1 in \cite{penbook}), $\Pr ( X < np/2 ) = O (\exp (-Cnp))$
for all $n, p$.  $\qed$

\begin{lemma}
Suppose $d \in \N$ and $\alpha > d/2$. For $\eps>0$ sufficiently small, 
we have that
\begin{align}
\label{sqa}
 \lim_{n \to \infty} \sup_{m: |n-m| \leq n^{(1/2)+\eps} } \left| \Exp [ \O^{d,\alpha} (\U_n)] - \Exp[ \O^{d,\alpha} (\U_m)] \right|
=0.\end{align}
\end{lemma}
\proof For ease of notation, let 
$\mu_n := \Exp[\O^{d,\alpha} (\U_n)]$.
By monotonicity of $\mu_n$,
\[ \sup_{m: |n-m| \leq n^{(1/2)+\eps} } | \mu_n - \mu_m | \leq \max \{ \mu_{n + \lceil n^{(1/2)+\eps} \rceil}
-\mu_n , \mu_n - \mu_{n - \lceil n^{(1/2)+\eps} \rceil} \},\]
so it suffices to show that both terms in the maximum tend to zero
as $n \to \infty$.
Consider the $\alpha \in (d/2,d)$ 
case of (\ref{ggjj}).
Now by Lemma \ref{xp}
we have, for small enough $\eps>0$,
\[ \sum_{i = n +1}^{n+\lceil n^{(1/2)+\eps} \rceil} h(i) 
\leq C n^{(1/2)+\eps} \sup_{i:n \leq i \leq n + \lceil n^{(1/2)+\eps} \rceil } h(i)
\leq C n^{(1/2)-(\alpha/d)-(1/d)+2\eps }
= o (n^{(1/2)-(\alpha/d)})  ,\]
which tends to $0$ 
 as $n \to \infty$,
given that $\alpha>d/2$.
Thus by (\ref{ggjj}), as $n\to\infty$,
\[ \mu_{n + \lceil n^{(1/2)+\eps} \rceil}
-\mu_n
= C n^{1-(\alpha/d)} 
( (1+n^{-(1/2)+\eps})^{1-(\alpha/d)} -1 ) + o(1),\]
for some $C \in (0,\infty)$. But this is $O(n^{(1/2)-(\alpha/d)+\eps})$, which tends
to zero for $\alpha >d/2$ and $\eps$ small enough. Similarly for
$\mu_n - \mu_{n - \lceil n^{(1/2)+\eps} \rceil}$. Thus we obtain (\ref{sqa}) for 
$\alpha \in (d/2,d)$.

Now suppose that $\alpha =d$. This time we have 
(\ref{ppiiuu});
by Lemma \ref{xp} the sum in (\ref{ppiiuu}) tends to $0$ as $m, \ell \to \infty$. 
Thus for $\eps>0$ small enough
\[ \mu_{n + \lceil n^{(1/2)+\eps} \rceil}
-\mu_n
= v_d^{-1} \log \left( \frac{ n + \lceil n^{(1/2)+\eps} \rceil}{n}
\right) +o(1) = O ( n^{\eps-(1/2)}) +o(1) \to 0,\] 
and similarly for
$\mu_n - \mu_{n - \lceil n^{(1/2)+\eps} \rceil}$.
 Thus we get (\ref{sqa}) for 
$\alpha =d$. The case $\alpha >d$ is straightforward, since there (see Proposition \ref{proprop})
$\mu_n \to \mu (d,\alpha) \in (0,\infty)$ as $n \to \infty$.
$\qed$\\

To prepare for the proof of Theorem \ref{onngthm},
we modify
 the technique used in the proof of Lemma \ref{mmmmm}
above.
For $n, m \in \N$ with $m <n$ set $Y^{(m)}_n :=
\tO^{d,\alpha} (\U_n) - \tO^{d,\alpha} (\U_m)$, i.e.~$Y^{(m)}_n$ is the centred
total weight of
edges in the ONG on $\U_n$ counting only edges from points after the first $m$ in the sequence. 
With $\F_i$ 
the $\sigma$-field generated by $(\bU_1,\ldots,\bU_i)$,
set
\[ D_i^{(n,m)} := \Exp [ Y^{(m)}_n \mid \F_i ] - \Exp [ Y^{(m)}_n \mid \F_{i-1} ],\]
so that for fixed $n,m$ the $D_i^{(n,m)}$ are martingale differences and
\[ Y^{(m)}_n = \sum_{i=1}^n D_i^{(n,m)} .\]
As in Section \ref{vars}, for $i \in \N$
let $\bU_i'$ be an independent copy of $\bU_i$. For  $i \leq n$ let
$\U_n^i$ be the sequence $\U_n$ but with 
$\bU_i$ 
replaced by $\bU_i'$. If $i > n$, we take $\U_n^i=\U_n$.  Define
\begin{align*}
 \Delta_i^{(n,m)} := [ \tO^{d,\alpha} (\U_n^i) 
-\tO^{d,\alpha} (\U_m^i) ] - [
\tO^{d,\alpha} (\U_n)-\tO^{d,\alpha} (\U_m) ]\\
= [\O^{d,\alpha} (\U_n^i) 
-\O^{d,\alpha} (\U_m^i) ]
- [ \O^{d,\alpha} (\U_n) - \O^{d,\alpha} (\U_m) ].\end{align*}
Then, similarly to before, 
\[  D_i^{(n,m)} = - \Exp [ \Delta_i^{(n,m)} \mid \F_i ].\]
Analogously to before, we decompose $\Delta_i^{(n,m)}$ into six parts.
For $i >m$,
let $\Delta_{i,1}^{(n,m)}$ be the weight of the edge from $\bU_i$,
and $\Delta_{i,2}^{(n,m)}$ be the weight of the edge from $\bU_i'$.
For $i \leq m$, set $\Delta_{i,1}^{(n,m)}=\Delta_{i,2}^{(n,m)}=0$.
For all $i$, let $\Delta_{i,\ell}^{(n,m)}$ for $\ell=3,4$
be the total weight of edges
incident to $\bU_i$, $\bU_i'$ respectively from $\{ \bU_{m+1}, \bU_{m+2}, \ldots, \bU_n \}$.
Let $\Delta^{(n,m)}_{i,5}$ be the total weight of edges in the
ONG on $( \bU_1, \ldots, \bU_{i-1}, \bU_{i+1}, \ldots, \bU_n)$ from points
in $\{\bU_{m+1},\ldots,\bU_n\}$ 
that are joined to $\bU'_i$ in the ONG on $\U^i_n$.
Let $\Delta^{(n,m)}_{i,6}$ be the total weight of edges in the
ONG on $( \bU_1, \ldots, \bU_{i-1}, \bU_{i+1}, \ldots, \bU_n)$ from points
in $\{\bU_{m+1},\ldots,\bU_n\}$ that are
 joined to $\bU_i$ in the ONG on $\U_n$.
Then we have 
\[ \Delta_{i}^{(n,m)} = \Delta_{i,2}^{(n,m)} + \Delta_{i,4}^{(n,m)}
+\Delta_{i,6}^{(n,m)}
 - \Delta_{i,1}^{(n,m)} 
-\Delta_{i,3}^{(n,m)} 
-\Delta_{i,5}^{(n,m)}.\]
Note that $\Delta_{i,\ell}^{(n,m)}
\geq \Delta_{i,\ell}^{(n,m+1)}$ and
$\Delta_{i,\ell}^{(n,1)}
= \Delta_{i,\ell}^{(n)}$ as defined in Section \ref{vars}.
Analogously 
to Lemma \ref{clem}
above, we have the following.

\begin{lemma}
\label{clem2}
For any $\alpha>0$
there exists $C \in (0,\infty)$ such that
for all $\ell \in \{1,\ldots,6\}$
\begin{align}
\label{claim1}
\Exp [ (\Exp [ \Delta_{i,\ell}^{(n,m)} \mid \F_i ])^2 ] \leq C i^{-2\alpha/d} ,\end{align}
for $m \leq i \leq n$, and, for $i < m \leq n$,
\begin{align}
\label{claim2}
\Exp [ (\Exp [ \Delta_{i,\ell}^{(n,m)} \mid \F_i ])^2 ] \leq C (\max \{m-i,i\})^{-2\alpha/d} .\end{align}
 \end{lemma}
\proof
The argument in Lemma \ref{clem} carries through, so that 
(\ref{claim1})  holds for all $i$. 
Indeed, $\Delta_{i,\ell}^{(n,m)}
\leq \Delta_{i,\ell}^{(n)}$ and so Lemma \ref{clem} implies (\ref{claim1})
for all $i,\ell$.
Thus to obtain (\ref{claim2})
we need to show that there exists $C \in (0,\infty)$
such that for all $\ell$ and all $i < m \leq n$
\begin{equation}
\label{aa1}
 \Exp[ (\Exp [ \Delta^{(n,m)}_{i,\ell} \mid \F_i ])^2 ] \leq 
C (m-i)^{-2\alpha/d}.\end{equation}

Thus suppose $i < m$. 
In this case, we need only consider
$\Delta_{i,\ell}^{(n,m)}$ for $\ell\geq 3$, since $\Delta_{i,\ell}^{(n,m)}=0$
for $\ell \in \{1,2\}$.
First take $\ell=3$, dealing with the edges incident to $\bU_i$.
There are $m-i$ points of $\U_n$
with mark (index) greater than $i$ but not more than $m$, and edges from
these points to $\bU_i$ are not counted in $\Delta_{i,3}^{(n,m)}$. 
Recall that $V_i, V_i'$ is the Voronoi cell of $\bU_i, \bU_i'$ respectively 
with respect to itself and
$\{ \bU_1,\ldots,\bU_{i-1}\}$, and $B_i$ 
is a minimal-volume 
 $d$-cube with $V_i \subseteq B_i \subseteq (0,1)^d$.

By an argument analogous to that in the proof
of Lemma \ref{clem},
discarding points of $\{\bU_{i+1},\ldots,\bU_m\}$
that fall outside $B_i$  can only increase $\Delta_{i,3}^{(n,m)}$.
It follows that, with the same notation as in that proof,
\begin{align*}
\Delta_{i,3}^{(n,m)} = \sum_{j=m+1}^n \1_{E(\bU_j,\bU_i; \U_j)}
 \| \bU_j - \bU_i \|^\alpha \leq \sum_{j: m+1 \leq j \leq n, \bU_j \in B_i}
  \1_{E(\bU_j,\bU_i; \U_{j,i})}
 \| \bU_j - \bU_i \|^\alpha .\end{align*}
Let
$M \sim {\rm Bin} (m-i, |B_i|)$ be the number
of points of $\{\bU_{i+1},\ldots,\bU_m\}$ that fall in $B_i$.
Thus
$\Delta^{(n,m)}_{i,3}$
is stochastically dominated by 
\[ \O^{d,\alpha}_{\bU_i} (B_i; n)-\O^{d,\alpha}_{\bU_i}(B_i;M) 
\eqd |B_i|^{\alpha/d} [ \O_\bx^{d,\alpha} 
(  n)
- \O^{d,\alpha}_\bx (   M)],\]
for some $\bx \in (0,1)^d$, by scaling.
Hence for some $C \in (0,\infty)$
\[ \Exp [ \Delta^{(n,m)}_{i,3} \mid \F_i ] \leq C |B_i|^{\alpha/d} 
\Exp[ (M+1)^{-\alpha/d} \mid \F_i  ],\]
by (\ref{ooky}). By Lemma \ref{binbin}, $\Exp [ (M+1)^{-\alpha/d} \mid \F_i ] 
 \leq C |B_i|^{-\alpha/d} (m-i)^{-\alpha/d}$,
so that for $i <  m$
\[ \Exp [ \Delta^{(n,m)}_{i,3} \mid \F_i ]
\leq  C (m-i)^{-\alpha/d} .\]
For $\ell=4$ a similar argument (with $V_i$ replaced by $V_i'$) holds.
Thus we obtain (\ref{aa1}) for $\ell \in \{3,4\}$.
For $\ell \in \{5,6\}$ a similar argument applies, using (\ref{ookyhat}) instead
of (\ref{ooky}) this time. $\qed$\\ 

\noindent
{\bf Proof of Theorem \ref{onngthm}.}
By Lemma \ref{clem2} and Cauchy--Schwarz we have 
that for $i < m$,
\[ \Exp [ (D_i^{(n,m)} )^2 ] \leq C (\max \{ m-i, i \})^{-2\alpha/d},\]
while for $i \geq m$, $\Exp [ (D_i^{(n,m)} )^2 ] \leq C  i ^{-2\alpha/d}$. 
Thus for $\alpha >0$, for $m < n$
\begin{align}
\label{mmcc}
 \Exp [ | \tO^{d,\alpha} (\U_n) - \tO^{d,\alpha} (\U_m) |^2 ] = \Exp [ (Y^{(m)}_n)^2 ]=
\sum_{i=1}^n \Exp[(D_i^{(n,m)})^2] \nonumber\\
\leq
C \sum_{i=1}^{\lceil m/2 \rceil} (m-i)^{-2\alpha/d} + C\sum_{i=\lfloor m/2 \rfloor}^m i^{-2\alpha/d}
+C \sum_{i=m+1}^n i^{-2\alpha/d}.\end{align}
In particular, for $\alpha >d/2$
the right-hand side of (\ref{mmcc}) is bounded by a constant times
$m^{1-(2\alpha/d)}$,
which tends to $0$ as $n$, $m$ tend to infinity.
Thus for $\alpha >d/2$,
$\tO^{d,\alpha}(\U_n)$ is  a Cauchy sequence in $L^2$, and hence as $n \to \infty$
it converges in $L^2$ to some limit random variable 
$Q(d,\alpha)$, with $\Exp[ Q(d,\alpha)]=\lim_{n \to \infty} \Exp[\tO^{d,\alpha}(\U_n)] =0$.  Thus we obtain (\ref{rrr}). 

Finally,
we prove the Poisson part (\ref{333}). 
As before, let 
$X_n := \O^{d,\alpha} (\U_n)$ and $\mu_n := \Exp[X_n]$. For $N(n)$ Poisson with mean
$n>0$, $\O^{d,\alpha} (\Po_n)$ has the distribution
of $X_{N(n)}$ and expectation $\Exp[ \mu_{N(n)} ]=:a_n$. 
Consider, for $n>0$
\begin{align}
\label{mmcm} \Exp [ | (X_{N(n)}-\mu_{N(n)}) - Q(d,\alpha) |^2 ]
 \leq  \sup_{m \geq n/2} \Exp [ |(X_m-\mu_m) - Q(d,\alpha) |^2 ] \nonumber\\ + 
\Exp [ |(X_{N(n)}-\mu_{N(n)} )-Q(d,\alpha)|^2 \1_{\{ N(n) < n/2\}}]
.\end{align}
For $\alpha >d/2$,
the $L^2$ convergence of $X_n-\mu_n$ to $Q(d,\alpha)$  (from (\ref{rrr}))
implies that the first term on the right-hand side
of (\ref{mmcm}) tends to zero, and that
the second term is bounded
by a constant times
$\Pr ( N(n) < n/2 )$, 
which tends to zero as $n \to \infty$. So, for $\alpha >d/2$, 
\begin{align}
\label{uuu}
 X_{N(n)} - \mu_{N(n)} \inLL Q(d,\alpha), ~{\rm as} ~ n \to \infty.\end{align}

First suppose $\alpha>d$. Here (see 
Proposition \ref{proprop})
$\mu_n \to \mu := \mu (d,\alpha)\in (0,\infty)$ as $n \to \infty$. 
It follows, by a similar argument to (\ref{mmcm}),
that $\mu_{N(n)}$ converges to $\mu$ in $L^2$ and
$a_n = \Exp[\mu_{N(n)}] \to \mu$ also. Thus, with (\ref{uuu}), as $n \to \infty$
\[ \tO^{d,\alpha} (\Po_n)
=X_{N(n)} - a_n = (X_{N(n)} - \mu_{N(n)} ) + (\mu_{N(n)} - \mu) + (\mu -a_n) \inLL Q(d,\alpha).\]

For $\alpha \in (d/2, d]$, $\mu_n \to \infty$ as $n \to \infty$.
Recall (\ref{iiii}).
With $p_m(n)=\Pr(N(n)=m)$,
the 
middle bracket
in (\ref{iiii}) satisfies,
for $\eps>0$,
\begin{align}
\label{tttz}
|a_n - \mu_{\lfloor n \rfloor}|
 = \sum_{m \in \N:|m-n|<n^{(1/2)+\eps}} |\mu_m-\mu_{\lfloor n \rfloor}| p_m(n)
+ \sum_{m \in \N:|m-n| \geq n^{(1/2)+\eps}} |\mu_m-\mu_{\lfloor n \rfloor}| p_m(n).\end{align}
Using the trivial bound $\mu_m \leq Cm$,  the second sum in (\ref{tttz}) 
is bounded by a constant times
\[ \sum_{m \in \N:|m-n| \geq n^{(1/2)+\eps}} (m+n) p_m(n) \leq  \Exp [ (N(n)+n) \1_{ \{ |N(n)-n|
 \geq n^{(1/2)+\eps} \}} ],
\]
which by Cauchy--Schwarz is bounded by 
\[ (\Exp [ (N(n)+n)^2])^{1/2} ( \Pr( |N(n)-n| \geq n^{(1/2)+\eps} ))^{1/2} \to 0, \]
as $n \to \infty$, by standard Poisson tail bounds.
The first sum in (\ref{tttz}) satisfies
\[
\sum_{m \in \N:|m-n|<n^{(1/2)+\eps}} |\mu_m-\mu_{\lfloor n \rfloor}| p_m(n)
\leq \sup_{m \in \N:|m-n|<n^{(1/2)+\eps}} |\mu_m-\mu_{\lfloor n \rfloor}|
,\]
which tends to zero as $n \to \infty$ by (\ref{sqa}). Thus
for $\alpha \in (d/2,d]$, as $n \to \infty$,
\begin{align}
\label{iiis}
|a_n - \mu_{\lfloor n \rfloor}| \to 0.\end{align} 
Also, from Lemma \ref{nevv} we have that, 
for $\alpha \in (d/2,d]$, $\Exp [ |\mu_{N(n)} - \mu_{\lfloor n \rfloor} |^2 ] \to 0$,
so that
 \begin{align}
 \label{i}
 \mu_{N(n)} - \mu_{\lfloor n \rfloor} \inLL 0, ~{\rm as} ~ n \to \infty.\end{align}
 Thus from (\ref{iiii}) with (\ref{uuu}), (\ref{iiis}) and (\ref{i}) we obtain
 the result for $\alpha \in (d/2,d]$ also.
$\qed$\\

\begin{center}{\bf Acknowledgements}
\end{center}

The author was partially supported
by the Heilbronn Institute for Mathematical Research.
Some of this work was carried out at the University of Bath. The author is
 grateful to Mathew Penrose
for many
helpful discussions and suggestions on the subject of this paper,
and also to an anonymous referee for a careful
reading of a previous version of this paper
and comments that have led to several improvements.


\begin{thebibliography}{99}

\bibitem{as} Abramowitz, M.~and Stegun, I.A.~(Eds.) (1965)
Handbook of Mathematical Functions, National Bureau of Standards,
Applied Mathematics Series {\bf 55}, U.S. Government Printing
Office, Washington D.C.

\bibitem{avbert} Avram, F.~and Bertsimas, D. (1993)
 On central limit
theorems in geometrical probability, \emph{Ann. Appl. Probab.}
\textbf{3} 1033--1046.

\bibitem{by} Baryshnikov, Yu.~and Yukich, J.E. (2005)
Gaussian limits
for random measures in geometric
probability, {\em Ann. Appl. Probab.} {\bf 15} 213--253.

\bibitem{bbcr} Berger, N., Bollob\'as, B., Borgs, C., Chayes,  J., and Riordan, O. (2003)
 Degree distribution
of the FKP network model. In {\em
Automata, Languages and Programming},
Lecture Notes in Computer
Science {\bf 2719},
Springer, Berlin, pp.~725--738.

\bibitem{bol1}
Bollob\'as, B.~and Riordan, O.M. (2003)
 Mathematical results on scale-free random graphs. In
  {\em Handbook of Graphs and Networks},  Wiley-VCH, Weinheim, pp.~1--34.

\bibitem{dorog}
Dorogovstev, S.N.~and Mendes, J.F.F. (2003)
Evolution of Networks,
 Oxford University Press, Oxford.

\bibitem{fkp} Fabrikant, A., Koutsoupias, E.~and Papadimitriou, C.H. (2002) 
Heuristically optimized trade-offs: a new paradigm for power laws in the internet. In
{\em Automata, Languages and Programming},
Lecture Notes in Computer Science {\bf 2380}, 
Springer, Berlin, pp.~110--122.  

\bibitem{huang} Huang, K. (1987)
Statistical Mechanics, 2nd ed., Wiley, New York.

\bibitem{KL} Kesten, H.~and Lee, S. (1996)
 The central limit
theorem for weighted minimal spanning trees on random points, {\em
Ann. Appl. Probab.} {\bf 6} 495--527.

\bibitem{neinrusch} Neininger, R.~and R\"uschendorf, L. (2004)
A general limit theorem for recursive algorithms and combinatorial structures,
{\em Ann. Appl. Probab.} {\bf 14} 378--418.

\bibitem{newman} Newman, M.E.J. (2003) 
The structure and function
of complex networks, {\em SIAM Rev.}  {\bf 45} 167--256.

\bibitem{penbook} Penrose, M. (2003) 
Random Geometric Graphs,
Oxford Studies in Probability {\bf 6}, Clarendon Press,
Oxford.

\bibitem{mdp} Penrose, M.D. (2005)
Multivariate spatial central limit theorems with applications
to percolation and spatial graphs, {\em Ann. Probab.} {\bf 33} 1945--1991.

\bibitem{mpgauss} Penrose, M.D. (2007)
Gaussian limits for random geometric measures,
{\em Electron. J. Probab.} {\bf 12} 989--1035.

\bibitem{total} Penrose, M.D.~and Wade, A.R. (2006)
 On the total
length of the random
 minimal directed spanning tree, {\em Adv. Appl. Probab.} {\bf 38} 336--372.

\bibitem{ong} 
Penrose, M.D.~and Wade, A.R. (2008)
Limit theory for the random on-line nearest-neighbor graph,
{\em Random Structures Algorithms} {\bf 32}
125--156.

\bibitem{py1} Penrose, M.D.~and Yukich, J.E. (2001)
 Central limit
theorems for some graphs in computational geometry, {\em Ann.
Appl. Probab.} {\bf 11} 1005--1041.


\bibitem{py2} Penrose, M.D.~and Yukich, J.E. 
(2003) Weak laws of large numbers in geometric probability,
{\em Ann. Appl. Probab.} {\bf 13} 277--303.

\bibitem{llns} Wade, A.R. (2007) 
Explicit laws of large numbers for random 
nearest-neighbour-type graphs, {\em Adv. Appl. Probab.}
{\bf 39} 326--342.
 
\bibitem{will} Williams, D. (1991)
Probability with Martingales,
Cambridge University Press, Cambridge.


\end{thebibliography}
\end{document}